\documentclass{amsart}

\usepackage{amssymb,amsfonts, amsmath, amsthm}
\usepackage[all,arc]{xy}
\usepackage{enumerate}
\usepackage{mathrsfs}
\usepackage{mathtools}
\usepackage[mathscr]{euscript}
\usepackage{array}
\usepackage{appendix}

\usepackage{tikz}
\usepackage{tikz-cd}
\usepackage{mathpazo}
\usepackage{textcomp}
\usepackage{bbold}
\usepackage[bbgreekl]{mathbbol}

\usepackage[pagebackref, colorlinks=true, linkcolor=blue, citecolor = red]{hyperref}
\renewcommand*{\backref}[1]{}
\renewcommand*{\backrefalt}[4]{({%
		\ifcase #1 Not cited.%
		\or On p.~#2%
		\else On pp.~#2%
		\fi%
	})}

\usepackage[nameinlink,capitalise,noabbrev]{cleveref}
\crefname{subsection}{Subsection}{Subsection}

\usetikzlibrary{patterns}

\usepackage{geometry}
\usepackage{todonotes}

\DeclareMathAlphabet{\mathbbe}{U}{bbold}{m}{n}

\def\DDelta{{\mathbbe{\Delta}}}
\newcommand{\DD}{\DDelta}

\newcommand{\A}{\mathscr{A}}
\newcommand{\C}{\mathscr{C}}
\newcommand{\D}{\mathscr{D}}
\newcommand{\E}{\mathscr{E}}
\newcommand{\Y}{\mathscr{Y}}
\newcommand{\G}{\mathscr{G}}
\newcommand{\sL}{\mathscr{L}}
\newcommand{\N}{\mathscr{N}}
\newcommand{\bN}{\mathbb{N}}
\newcommand{\s}{\mathscr{S}}

\newcommand{\cM}{\mathcal{M}}
\renewcommand{\O}{\mathscr{O}}
\newcommand{\R}{\mathscr{R}}
\newcommand{\K}{\mathscr{K}}
\newcommand{\cO}{\mathcal{O}}
\newcommand{\Q}{\mathscr{Q}}
\newcommand{\U}{\mathscr{U}}
\newcommand{\V}{\mathscr{V}}

\newcommand{\sP}{\mathscr{P}}
\newcommand{\cP}{\mathcal{P}}
\newcommand{\X}{\mathscr{X}}
\renewcommand{\Y}{\mathscr{Y}}

\newcommand{\Map}{\mathrm{Map}}
\newcommand{\Hom}{\mathrm{Hom}}

\newcommand{\Sub}{\mathrm{Sub}}

\newcommand{\id}{\mathrm{id}}
\newcommand{\comma}{,}

\newcommand{\idem}{\mathrm{idem}}
\newcommand{\Prop}{\mathrm{Prop}}
\newcommand{\Idem}{\mathrm{Idem}}

\newcommand{\Fun}{\mathrm{Fun}}

\newcommand{\Fib}{\mathscr{F}\mathrm{ib}}

\newcommand{\colim}{\mathrm{colim}}

\newcommand{\Und}{\mathscr{U}\mathrm{nd}}
\newcommand{\Lift}{\mathscr{L}\mathrm{ift}}
\newcommand{\LCCC}{\mathrm{LCCC}}

\newcommand{\RFib}{\mathcal{R}\mathcal{F}\mathrm{ib}}
\newcommand{\Cart}{\mathcal{C}\mathrm{art}}

\newcommand{\otherComp}{\mathscr{C}\mathrm{omp}}
\newcommand{\SegCat}{\mathscr{S}\mathrm{eg}\mathscr{C}\mathrm{at}}

\newcommand{\Set}{\mathscr{S}\mathrm{et}}
\newcommand{\cat}{\mathscr{C}\mathrm{at}}

\newcommand{\ds}{\displaystyle}

\newcommand{\CSS}{\mathscr{C}\mathscr{S}\mathscr{S}}
\newcommand{\CSU}{\mathscr{C}\mathscr{S}\mathscr{U}}
\newcommand{\sO}{\mathscr{O}}
\newcommand{\Univ}{\mathscr{U}\mathrm{niv}}

\newcommand{\CSUniv}{\mathscr{C}\mathscr{S}\mathscr{U}\mathrm{niv}}
\newcommand{\all}{\mathrm{all}}
\newcommand{\Oall}{\sO^{(\all)}}
\newcommand{\Aall}{\A^{(\all)}}
\newcommand{\OS}{\sO^{(S)}}

\newcommand{\Acc}{\mathscr{A}\mathrm{cc}}
\newcommand{\Log}{\mathscr{L}\mathrm{og}}
\newcommand{\Ind}{\mathrm{Ind}}
\newcommand{\cmp}{\mathrm{cmp}}
\newcommand{\sma}{\mathrm{sma}}

\newcommand{\Accsma}{\Acc^{\sma}}
\newcommand{\Elem}{\mathscr{E}\mathrm{lem}}

\newcommand{\AccElem}{\Acc\Elem}

\newcommand{\pbsq}[8]{
  \begin{tikzcd}[row sep=0.3in, column sep=0.3in]
    #1 \arrow[r, "#5"] \arrow[d, "#6"'] \arrow[dr, phantom, "\ulcorner", very near start]
    \pgfmatrixnextcell #2 \arrow[d, "#7"] \\
    #3 \arrow[r, "#8"']
    \pgfmatrixnextcell #4
  \end{tikzcd}
}

\newcommand{\fibref}[1]{\cref{sec:infty cats}(\ref{#1})} 
\newcommand{\citenno}[1]{\cite[#1]{rasekh2021nno}}
\newcommand{\citetrunc}[1]{\cite[#1]{rasekh2018truncations}}

\newtheorem{theone}[equation]{Theorem}
\newtheorem*{thetwo}{Theorem}
\newtheorem{lemone}[equation]{Lemma}
\newtheorem{propone}[equation]{Proposition}
\newtheorem{corone}[equation]{Corollary}

\theoremstyle{definition}
\newtheorem{defone}[equation]{Definition}
\newtheorem{exone}[equation]{Example}

\newtheorem*{nottwo}{Note on Terminology}

\theoremstyle{remark}
\newtheorem{remone}[equation]{Remark}

\newtheorem{notone}[equation]{Notation}

\newtheorem{conjone}[equation]{Conjecture}
\newtheorem{warnone}[equation]{Warning}

\numberwithin{equation}{section}

\newtheoremstyle{TheoremNum}
{}{}              % space between body and thm
{\itshape}                      % Thm body font
{}                              % Indent amount (empty = no indent)
{\bfseries}                     % Thm head font
{.}                             % Punctuation after thm head
{ }                             % Space after thm head
{\thmname{#1}\thmnote{ \bfseries #3}}% Thm head spec
\theoremstyle{TheoremNum}
\newtheorem{thmn}{Theorem}

\makeatletter
\def\@seccntformat#1{%
  \expandafter\ifx\csname c@#1\endcsname\c@section\else
  \csname the#1\endcsname\quad
  \fi}
\makeatother

\title{A Theory of Elementary Higher Toposes}

\author{Nima Rasekh}

\date{January 2022}

\address{{\'E}cole Polytechnique F{\'e}d{\'e}rale de Lausanne, SV BMI UPHESS, Station 8, CH-1015 Lausanne, Switzerland}

\email{nima.rasekh@epfl.ch}

\begin{document}

\begin{abstract}
We define an elementary $\infty$-topos that simultaneously generalizes an elementary topos and Grothendieck $\infty$-topos. We then prove it satisfies the expected topos theoretic properties, such as descent, local Cartesian closure, locality and classification of univalent morphisms, generalizing results by Lurie \cite{lurie2009htt} and Gepner-Kock \cite{gepnerkock2017univalence}. We also define $\infty$-logical functors and show the resulting $\infty$-category is closed under limits and filtered colimits, generalizing the analogous result for elementary toposes and Grothendieck $\infty$-toposes. Moreover, we give an alternative characterization of elementary $\infty$-toposes and their $\infty$-logical functors via their ind-completions. Finally we generalize these results by discussing the case of elementary $(n,1)$-toposes and give various examples and non-examples.
\end{abstract}

\maketitle
\addtocontents{toc}{\protect\setcounter{tocdepth}{1}}

\tableofcontents

 \section{Introduction}
 
 \subsection{Set Theory as a Foundation for Mathematics}
 Since the work of Cantor in the late 19th century \cite{cantor1874settheory}, set theory has established itself as a standard foundation for broad swaths of modern mathematics. As a result it is now fairly standard that in order to define a mathematical object, such as group, variety, manifold, or ..., one starts with a set and then adds the desired list of conditions. Having this common foundation allows mathematicians all over the world to effectively communicate ideas and introduce new objects.
 
 The fact that set theory can now be comfortably used to do mathematics follows a century long development using various mathematical foundations.  One such approach is axiomatic set theory first developed by Frege \cite{frege1893arithmetik}, Zermelo \cite{zermelo1908settheory} and Fraenkel \cite{fraenkel1922settheory}, leading to the now famous {\it Zermelo-Fraenkel} (ZF) set theory. Another approach arose in the branch of categorical logic and originated with work of Lawvere \cite{lawvere1964elementarysets} and Tierney \cite{tierney1973elementarytopos}, leading to a foundation via {\it elementary toposes}. Another, third approach, comes from {\it type theory}, which has arose in the work of Russel \cite{russell1908typetheory}, as a way to avoid paradoxes in early versions of axiomatic set theory and is particularly suitable in constructive and computational mathematics \cite{church1940typetheory,martinlof1975inttypetheoriespredicative}. Moreover, there are ways to translate between these foundations: Every (higher-order) type theory gives us an elementary topos and vice versa \cite{lambekscott1988higherorderlogic}. Moreover, we can construct set theories both out of elementary topos theories \cite[Subsection VI.10]{maclanemoerdijk1994topos} and type theories \cite{aczel1978setvstype}.

The key aspect of set theory as a foundation for mathematics is not what a set {\bf is}, but rather how it {\bf behaves} in the sense that we want to know which axioms our chosen foundation satisfies. The development above allows mathematicians of all stripes to choose the axioms that fits their needs to develop their mathematical theory, which might or might not include well-known axioms such as the {\it axiom of choice}, the {\it continuum hypothesis}, the {\it law of excluded middle}, and many other axioms, being assured of soundness of the theory without having to necessary grapple with the underlying foundation and which one of the three is used to construct a specific model.
 
 \subsection{Homotopical Mathematics}
 Homotopical mathematics has its starting point in the work of Poincar{\'e}, who recognized that we can distinguish two topological spaces by analyzing homotopy classes of based loops, introducing the first homotopy invariant \cite{poincare1895analysissitus}. This led to further invariants, such as {\it higher homotopy groups} \cite{cech1932homotopygroups} and {\it homology} and {\it cohomology theories} \cite{eilenbergsteenrod1945axioms}. 
 
 The development of these invariants motivated the introduction of up-to-homotopy structures, such {\it $A_\infty$-spaces}, which are monoids in spaces where the composition and associativity only hold up to a choice of coherent homotopies, a prominent example being the based loop space \cite{stasheff1963ainftyspace}. This example and many other illustrate why set theoretical foundations are not suitable for the study of homotopical mathematics, as we would rather like to use a framework in which the loop space is directly a monoid (in the sense that the structure of a monoid is itself only defined up to coherent homotopies). This vision spurred a journey, which led to the development of various notions of {\it homotopical} or {\it weak categories}, such as {\it quasi-categories} \cite{boardmanvogt1973qcats}, {\it complete Segal spaces} \cite{rezk2001css}, {\it Kan enriched categories} \cite{bergner2007bergnermodelcat} and others, which are all now considered examples of the theory of {\it $(\infty,1)$-categories} (or simply {\it $\infty$-categories}) \cite{bergner2010survey}, the standard setting for homotopy coherent mathematics. 
  
 Similar to the previous subsection homotopical mathematical requires its own basic building block that can be used to characterize different objects of interest and we could analogously call a {\it theory of spaces}\footnote{To avoid ambiguity there have also been suggestions to call it {\it animated set theory} \cite{cesnaviciusscholze2019flatcohomology}}. Again the key value is to determine how a space behaves rather than what a space is. Unlike the set-theoretical setting this foundation has not been fully established yet, however, we are currently witnessing its development. Ideally we would expect that the three foundations outlines above, namely axiomatic set theory, elementary topos theory and type theory, have a homotopy invariant analogues that also interact with each other in a similar manner.
  
 A first step important step towards realizing this vision is a homotopy invariant analogue of type theory, known as {\it homotopy type theory}, which started with the work of  Awodey, Warren \cite{awodeywarren2009identitytypes} and Voevodsky \cite{voevodsky2014origins}, and now covers many fundamental aspects of topology, such as fundamental groups, Blakers-Massey and truncations \cite{hottbook2013}.
 
\subsection{Towards Elementary \texorpdfstring{$\infty$}{oo}-Toposes}
 The goal of this paper is to establish another pillar of homotopy invariant foundation, by generalizing elementary toposes to an appropriate notion of {\it elementary $\infty$-topos}. In order to understand how we can obtain a working definition, it is instructive to review how the definition of elementary toposes arose in the first place. 
 
 The story of topos theory starts with the work of Grothendieck and more generally Bourbaki, who needed a proper categorical framework to study algebro-geometric objects, introducing {\it Grothendieck topologies} and {\it sheaves} and hence defining a {\it Grothendieck topos} as a category of sheaves \cite{sga1972tome1}. Lawvere and Tierney recognized that Grothendieck toposes can be an appropriate framework for categorical logic and an axiomatic study of sets if they could avoid the set-theoretical assumptions inherent in the definition of Grothendieck toposes (a technical condition known as {\it local presentability}). Their solution was to add a certain universal object, the {\it subobject classifier}, which allowed them to recover the desired results without assuming local presentability  \cite{lawvere1964elementarysets,tierney1973elementarytopos}, and defining elementary toposes as {\it Cartesian closed categories} with subobject classifier and now many important results can now be found in the famous work of Johnstone \cite{johnstone2002elephanti, johnstone2002elephantsii}.
 
 The $(\infty,1)$-categorical story of topos theory starts again very similarly, namely with the desire by Lurie to study objects in {\it derived algebraic geometry}, motivating him to define {\it Grothendieck $\infty$-toposes} as appropriate generalizations of Grothendieck toposes \cite{lurie2009htt}. Around the same time, Rezk also developed the equivalent theory of {\it model toposes}, with the intention of better understanding the concept of {\it descent}, in particular introducing the (now standard) idea that Grothendieck $(\infty,1)$-toposes are locally Cartesian closed presentable $(\infty,1)$-categories that satisfy descent \cite{rezk2010toposes}. The historical analogy would suggest that we can obtain our desired generalization by replacing the presentability condition with appropriately chosen universal objects. 
 
 Before we explain how to choose the universal object, we will make a short detour into various characterizations of categories and $(\infty,1)$-categories. The data of a category can be characterized in two equivalent ways: Either as a set of objects along with a hom set for any two objects that has a composition, or as a set of objects and set of morphisms along with functions corresponding to source, target, composition, ... \cite[Subsection 1.1]{riehl2017context}. We can think of the first as a {\it local} characterization of categories and the second as a {\it global} characterization. 
 
 This philosophy carries over to the $(\infty,1)$-categorical setting. There we have models of $(\infty,1)$-categories that follow the local philosophy, such as relative categories \cite{barwickkan2012relativecategory}, Kan enriched categories \cite{bergner2007bergnermodelcat} or topologically enriched categories \cite{lurie2009htt}, and models that follow the global philosophy, such as quasi-categories \cite{boardmanvogt1973qcats}, complete Segal spaces \cite{rezk2001css} and Segal categories \cite{bergner2007threemodels}. There are prominent ways of translating between the local and global perspective, such as $(\mathfrak{C}, N)$ between quasi-categories and Kan enriched categories \cite[Theorem 2.2.5.1]{lurie2009htt}, or $(F,R)$ between Segal categories and Kan enriched categories \cite[Theorem 8.6]{bergner2007threemodels}, $(K_\xi,N_\xi)$ between complete Segal spaces and relative categories \cite[Theorem 6.1]{barwickkan2012relativecategory}.

  \subsection{Elementary \texorpdfstring{$\infty$}{oo}-Topos: Local vs. Global}
 Having reviewed local and global approaches to $(\infty,1)$-categories, we claim that the appropriate notion of universal object is given by an internal $\infty$-category that classifies the morphisms over a given object. This universal object can be characterized using both the global perspective and the local perspective:
 \begin{itemize}
 	\item {\bf Local:} A {\it universe}, meaning an object $\U$, such that morphisms $ X \to \U$ correspond (homotopically uniquely) to a space of morphisms over $X$, along with local Cartesian closure. Here $\U$ should be thought of as the objects of our internal $\infty$-category and the local Cartesian closure gives us the morphism objects. 
 	\item {\bf Global:} A {\it complete Segal universe}, meaning an internal complete Segal object  $\U_\bullet$, such that morphisms $X \to \U_\bullet$ correspond (homotopically uniquely) to an $\infty$-category of objects over $X$. Here, we are using the fact that the theory of complete Segal spaces can be internalized to give us a theory of internal $\infty$-categories via complete Segal objects \cite{rasekh2017cartesian}.
 \end{itemize}
  
 We claim that similar to before the local and global approach do in fact coincide. Concretely, we can combine (and slightly simplify) \cref{def:eht} and \cref{prop:eht lccc}, fundamentally relying on \cref{the:main univ}, to give the following result.
 
 \begin{thetwo}
 	Let $\C$ be a finitely (co)complete $\infty$-category. Then the following coincide:
 	\begin{itemize}
 		\item It is locally Cartesian closed and has sufficient universes\footnotemark{}.
 		\item It has sufficient complete Segal universes\footnotemark[\value{footnote}]. \footnotetext{Here we are using several (complete Segal) universes in order to avoid size-related paradoxes.}
 	\end{itemize}
 An $\infty$-category that satisfies these equivalent conditions and has a subobject classifier is an elementary $\infty$-topos.
\end{thetwo}

 Using these characterizations we can deduce various valuable facts about elementary $\infty$-toposes. First of all we can study classify all univalent morphisms in an elementary $\infty$-topos, generalizing an analogous result by Gepner and Kock \cite[Corolary 3.11]{gepnerkock2017univalence}. In fact we can use univalent morphisms to give an alternative characterization of elementary $\infty$-toposes, which is closely related to the approach in \cite[Remark 3.27, Remark 4.21]{stenzel2020comprehension}. Moreover, using the local characterization, we realize that the definition (with some closure properties on the universe) coincides with the one suggested by Shulman \cite{shulman2017eht}.
 
 We can also witness that elementary $\infty$-toposes satisfy conditions known (or not known) from elementary topos theory. For example, similar to elementary toposes, we can prove the {\it fundamental theorem of $\infty$-topos theory}, which proves that the overcategory of an elementary $\infty$-topos is again an elementary $\infty$-topos (\cref{the:eht fundamental}). Moreover, using \cite{rasekh2021nno}, we can deduce that every elementary $\infty$-topos has a natural number object (\cref{the:eht nno}). Finally, we can also use \cite{rasekh2018truncations}, we can deduce that every elementary $\infty$-topos (with sufficiently closed universes), has a truncation functor (\cref{the:eht ntrunc}), generalizing the epi-mono factorization in every elementary topos. These results pave the way to tackle future problems, as we have discussed in greater detail in \cref{subsec:future}

 Having a definition of elementary $\infty$-toposes, we can also define an appropriate notion of functor, the {\it $\infty$-logical functor}. We can again use the local vs. global principle and combine \cref{def:infty logical functor} and \cref{prop:logical is lccc} to give the following result.
 
 \begin{thetwo}
 	Let $F:\C \to\D$ be a functor of elementary $\infty$-toposes that preserves finite limits and colimits. Then the following coincide:
 	\begin{itemize}
 		\item $F$ preserves the local Cartesian structure and universes.
 		\item $F$ preserves complete Segal universes.
 	\end{itemize}
 	A functor that satisfies these equivalent conditions and preserves the subobject classifier is an $\infty$-logical functor.
 \end{thetwo}
 This generalizes logical functors of elementary toposes  (\cref{def:et}) and results in the $\infty$-category $\Log_\infty$. It is known that various notions of toposes and their functors are closed under limits and filtered colimits. Indeed for elementary toposes this can be found in \cite[Page 218]{maclanemoerdijk1994topos}, and for Grothendieck $\infty$-toposes in \cite[Proposition 6.3.2.3, Theorem 6.3.3.1]{lurie2009htt}. We continue this tradition by proving the following result.
 
 \begin{thmn}[\ref{the:filtered colimit eht}]
 	The inclusion functor $\Log^{\cP}_\infty \to \widehat{\cat}_\infty$ preserves small limits and filtered colimits.
 \end{thmn}
 
 There is a stronger result due to Dubuc and Kelly, that proves that the category of elementary toposes and logical functors is also cocomplete (and in fact a category of algebras of a finitary monad over $\cat$) \cite{dubuckelly1983toposcomplete}. This in particular implies the existence of an initial object, which is known as the {\it free elementary topos} and plays a key role in the translation between type theories and elementary toposes \cite[Example D4.3.14]{johnstone2002elephantsii} \cite{lamekscott1980freetopos,lambek1980typessets,lambekscott1988higherorderlogic}. Given our desire to better understand this connection, we can ask ourselves whether there is an {\it free elementary $\infty$-topos}. While it is expected to be true, we shall postpone it to later and work and leave it as \cref{conj:free topos}.
 
 Even more generally, we can ask ourselves whether the inclusion functor $\Log_\infty \to \widehat{\cat}_\infty$ has a left adjoint, that takes every category to its {\it freely generated elementary $\infty$-topos}, which would also give us the free elementary $\infty$-topos if we start with the initial $\infty$-category. Again, we shall postpone it to later and work and leave it as \cref{conj:free topos adj}.
 
 While most of our work has focused on studying $\infty$-topos theory, there is also a relevant theory of $(n,1)$-toposes. The Grothendieck version was also developed by Lurie \cite[Section 6.4]{lurie2009htt}, and in \cref{sec:eht n} we introduce a theory of {\it elementary $(n,1)$-toposes} (\cref{def:eht n}) and their $(n,1)$-logical functors (\cref{def:logical n}), however only using the local approach (\cref{rem:eht not global}). 
 
 Despite the fact that an elementary $(n,1)$-topos is not an elementary $(m,1)$-topos (\cref{prop:eht m not n}), we can in fact obtain an elementary $(n+1,1)$-topos as the subcategory of $n$-truncated objects (\cref{the:truncation eht}). While it is known that every Grothendieck $(n,1)$-topos can be obtained this way \cite[Theorem 6.4.1.5]{lurie2009htt}, it is not known whether this holds for elementary $(n,1)$-toposes and we leave the relevant conjecture (\cref{conj:truncated eht}) for future work.
 
 \subsection{External Universes and Models of Homotopy Type Theory}
 Having established a theory of elementary $\infty$-toposes one would similarly expect that it corresponds to homotopy type theories, the same way that elementary toposes correspond to higher-order type theories. The main challenge towards realizing this goal is to show that the universes in an elementary $\infty$-topos are as strict as required by type theory, which is complicated by the fact that there are many different universes that need to be strictified. 
 
 In the setting of Grothendieck $\infty$-toposes, this issue has been solved by Shulman \cite{shulman2019inftytoposunivalent} by embedding the Grothendieck $\infty$-topos in a category of groupoid valued presheaves and constructing a single external universe, that can then be shown to satisfy the desired strictness results. This approach fundamentally hinges on the fact that Grothendieck $\infty$-toposes are all presentable and so can be modeled by model categories \cite{dugger2001universal}. Arbitrary elementary $\infty$-toposes will not satisfy such strictness conditions and so we cannot hope to recover such strong results directly. 
 
 Hence, as a first step we would like to find a way to characterize elementary $\infty$-toposes in a way that requires only a single universe. One possible motivation could come from \cite{abss2014etuniverse}, where the authors use algebraic set theory \cite{joyalmoerdjik1995algebraicset} to associate to each elementary elementary topos a {\it category of ideals} using the {\it directed system of inclusions}. We will leave the $\infty$-categorical approach to algebraic set theory to future work and instead use the ind-completion of an $\infty$-category to obtain the desired characterization with a single universe. In fact, in \cref{the:main external stuff} we prove a stronger result by giving a correspondence between various structured $\infty$-categories (such as universes, local Cartesian closure and subobject classifiers) and their ind-completions.
 
 One immediate benefit of this correspondence is that, despite their complicated appearances, the characterizations of elementary $\infty$-toposes and their $\infty$-logical functors via sufficient universes and their preservation are in fact quite sensible, as they do correspond to the existence of a single universe and its preservation (\cref{rem:justification}). Beyond that, the hope is that this alternative characterization can be employed in future work to prove that elementary $\infty$-toposes model homotopy type theories by strictifying the unique universe, generalizing the result by Shulman \cite{shulman2019inftytoposunivalent}.

 \subsection{Where to go from here} \label{subsec:future}
 While we address some questions regarding elementary $\infty$-topos theory, many questions remain.
 
 \begin{enumerate}
 	\item {\bf Free Algebras and $\infty$-Operads:} One elegant result in elementary topos theory is the existence of free finitary algebras, and in particular free monoids, using the notion of {\it $W$-types} \cite[Theorem D5.3.5]{johnstone2002elephantsii}. The construction fundamentally relies on the existence of natural number objects, which suggests the possibility of similar free constructions in elementary $\infty$-topos theory, however, the $\infty$-categorical nature raises several new challenges that need to be addressed.
 	
 	Unlike in the $1$-categorical case, in an $\infty$-category $\E$, making an object $E$ into a monoid involves constructing an infinite tower of morphisms $\mu_n: E^n \to E$, which are appropriately compatible. This data is often managed via an appropriate choice of $\infty$-operad and the monoid is then an algebra for such $\infty$-operads \cite{lurie2017ha}. Given that an $\infty$-operad incorporates the data of operations $\mu_n$ for all $n \in \bN$, the current definition relies on the natural numbers, which in an arbitrary elementary $\infty$-topos can differ, and, concretely, involve non-standard natural numbers \cite{rasekh2021filterquotient}. As a result, before we can even talk about constructing free algebras, the first important step is to give a definition of $\infty$-operads internal to an elementary $\infty$-topos $\E$ based on the internal natural numbers of $\E$.
 	\item {\bf Stabilizations:} One key development of algebraic topology is stable homotopy theory. It involves defining the $\infty$-category of spectra. Hence we would want a stabilization of an arbitrary elementary $\infty$-topos. There is an established notion of a stabilization of an $\infty$-category, however, in this case we face the same challenges as before, as the stabilization is indexed by the natural numbers \cite[Definition 1.4.2.8]{lurie2017ha}, which again can differ in an arbitrary elementary $\infty$-topos, and can in particular result in the existence of additional spheres. So what we would want is defining a stabilization operation that takes the additional natural numbers and spheres into account. 
 	\item {\bf Constructing Finite Colimits:} Another elegant result in elementary topos theory is the fact that we can construct finite colimits using Cartesian closure and the subobject classifier and so, unlike the $\infty$-categorical analogue, the definition of elementary toposes does not assume the existence of finite colimits \cite{mikkelsen1972finite,mikkelsen1976lattice,pare1974colimits}. One relevant question is whether a similar result holds in the $\infty$-categorical setting. As a first step, in joint work with Jonas Frey, we show that initial objects and finite coproducts can be constructed in every locally Cartesian closed $\infty$-category with subobject classifier, and so in particular in every elementary $\infty$-topos \cite{freyrasekh2021coprod}. The remaining step towards recovering all finite colimits is to construct coequalizers using universes.
 	\item {\bf Localizations:} Another fascinating result about elementary toposes is our ability to construct all left-exact localizations of an elementary topos via certain endomorphisms of the subobject classifier, the so-called {\it Lawvere-Tierney topologies} \cite[Section V]{maclanemoerdijk1994topos} and we would like to obtain a similar elegant characterization for elementary $\infty$-toposes. We do have some results in the case of Grothendieck $\infty$-toposes \cite{lurie2009htt,abfj2022topology,stenzel2021hottest}, however, obtaining an appropriate generalization is greatly complicated by the fact that even in Grothendieck $\infty$-toposes not every left-exact accessible localization can be obtained via Grothendieck topologies, as they can have $\infty$-connected morphisms, and the ones that can be obtained are known as {\it topological localizations} \cite[Definition 6.2.1.4]{lurie2009htt}. This suggests two reasonable next steps, that we should pursue:
 	\begin{itemize}
 		\item Proving that we can we can characterize topological localizations of elementary $\infty$-toposes via Lawvere-Tierney topologies on the subobject classifier.
 		\item Proving that more general left-exact localizations can be classified by certain choice of topology on a specific universe.
 	\end{itemize}
 	\item {\bf Models of Homotopy Type Theory:} Finally, as explained above, elementary toposes can be understood as models of various higher order type theories \cite[Proposition D.4.3.15]{johnstone2002elephantsii} and it is expected that there are similar connections between intensional type theories and $\infty$-categories, and particularly homotopy type theories with univalent universes and elementary $\infty$-toposes. 
 	
 	Making this connection precise is quite challenging as $\infty$-categories are not strict enough to model type theories directly and all existing results use the fact that we can strictify the $\infty$-category, either to model categories \cite{shulman2019inftytoposunivalent} or at least fibration categories \cite{kapulkinszumilo2019completequasicat}. Hence one first step towards realizing this goal is to give a definition of elementary $\infty$-toposes in a strict model of $\infty$-categories, such as possibly fibration categories. 
 \end{enumerate}
 
 \subsection{Background}
 We will assume general familiarity with the world of $\infty$-categories (any model) as presented in \cite{riehlverity2021elements,rezk2017qcats,lurie2009htt} and will only review some key concepts in \cref{sec:infty cats}. Some familiarity with $\infty$-topos theory as presented in \cite[Section 6]{lurie2009htt} or \cite{rezk2019leeds} would be helpful, however, we have reviewed most key results that we need in \cref{sec:review}.
 
 \subsection{Acknowledgment}
 I want to thank my advisor Charles Rezk for suggesting this topic and for his helpful comments, in particular regarding complete Segal universes. I also want to thank Mike Shulman for many fruitful conversations that contributed to \cref{the:main univ}. Moreover, I want to thank Valery Isaev for helpful comments regarding closed universes now explained in \cref{rem:comparison shulman}. I finally want to thank Nicola Gambino for pointing me to \cite{abss2014etuniverse}, which led to the material in \cref{sec:external universe}.

\section{From \texorpdfstring{$\infty$}{oo}-Categories to Representable Cartesian Fibrations} \label{sec:infty cats}
We will use the language of $\infty$-categories throughout. For that we will primarily rely on the theory of {\it $\infty$-cosmoses} introduced by Riehl and Verity \cite{riehlverity2017inftycosmos,riehlverity2021elements} and to some extent on the theory of {\it quasi-categories} as studied by Joyal \cite{joyal2008notes,joyal2008theory} and Lurie \cite{lurie2009htt,lurie2017ha}. Finally, some results depend on the theory of {\it representable Cartesian fibrations} introduced in  \cite{rasekh2017cartesian} and further analyzed in \cite{rasekh2021cartfibcss,rasekh2021cartfibmarkedvscso} and the associated study of {\it univalence} \cite{rasekh2021univalence}. We will review some of the key definitions and results here:

\begin{enumerate}
	\item \label{item:universes} We fix three universes, which we call {\it small}, {\it large} and {\it very large}.
	\item \label{item:cosmos} We denote the (very large) $\infty$-cosmos of (large) quasi-categories \cite{joyal2008notes,joyal2008theory} by $\Q\cat$. A {\it theory of $\infty$-categories} is an $\infty$-cosmos $\K$ such that the {\it underlying quasi-category functor} $(-)_0= \Hom_\K(1,-): \K \to \Q\cat$ is a biequivalence of $\infty$-cosmoi with inverse the tensor $1 \otimes -: \Q\cat \to \K$ (also called $\infty$-cosmos of $(\infty,1)$-categories \cite[Definition 1.3.10]{riehlverity2021elements}). According to \cite[Example 1.3.9]{riehlverity2021elements}, examples include the $\infty$-cosmos of {\it quasi-categories} itself, but also the $\infty$-cosmos of {\it complete Segal spaces} $\CSS$ \cite{rezk2001css,joyaltierney2007qcatvssegal}, {\it Segal categories} $\SegCat$ \cite{bergner2007threemodels}, and {\it $1$-complicial sets} $1-\otherComp$ \cite{lurie2009htt}. 
	\item \label{item:limits}
	We define finitely complete  $\infty$-categories as defined in \cite{riehlverity2021elements}. In particular, by \cite[Proposition 4.3.1]{riehlverity2021elements}, the limit of a diagram $F:I \to \C$ is given by the terminal object in the pullback $\infty$-category 
	$$1 \xrightarrow{ \ F \ } \C^I \xleftarrow{ \ \Delta \ } \C.$$  
   Note that by \cite[Appendix F]{riehlverity2021elements} in the particular case of quasi-categories the definition coincides with alternative definitions given in \cite{lurie2009htt,gepnerkock2017univalence}. Moreover, we use the notation ``{\it finitely (co)complete $\infty$-category}" to denote an $\infty$-category that has finite limits and colimits. Similarly, we use the notation ``{\it finitely (co)continuous}" to denote a functor that preserves finite limits and colimits.  
    \item \label{item:lccc} We say a finitely complete $\infty$-category $\C$ is locally Cartesian closed  if the pullback functor $f^*:\C_{/y} \to \C_{/x}$ has a right adjoint \cite[Definition 2.1.1]{riehlverity2021elements} for every morphism $f:x \to y$ in $\C$. Again, by \cite[Appendix F]{riehlverity2021elements} in the particular case of quasi-categories the definition coincides with the other common definitions \cite{lurie2009htt,gepnerkock2017univalence}. 
	\item \label{item:css} If $\K$ is a theory of $\infty$-categories, there is a biequivalence $\CSS=\mathrm{nerve} (\Hom_\K(1,-)):\K \to \CSS$ \cite[Example 1.3.9]{riehlverity2021elements}, that takes an $\infty$-category in $\C$ to its {\it underlying complete Segal space}.
	\item \label{item:inftycategories} We denote the underlying (very large) quasi-category of $\K$ by $\widehat{\cat_\infty}$ (which is not an object in the $\infty$-cosmos $\Q\cat$) and up to equivalence does not depend on the choice of $\infty$-cosmos. In particular, for two $\infty$-categories $\C,\D$ we denote the Kan complex of functors by $\Map_{\cat_\infty}(\C,\D)$. Moreover, we denote the full subcategory of (very large) $\infty$-groupoids (which we also call spaces) by $\widehat{\s}$. 
	\item \label{item:underlyinggroupoid} The inclusion functor $\widehat{\s} \hookrightarrow \widehat{\cat_\infty}$ has a right adjoint, which takes an $\infty$-category to its {\it underlying maximal $\infty$-groupoid}, which we denote by $(-)^\simeq: \widehat{\cat_\infty} \to \widehat{\s}$.
	\item \label{item:truncated} Recall a morphism $f:c \to d$ in a finitely complete $\infty$-category $\C$ is {\it $(-2)$-truncated} if it is an equivalence and {\it $n$-truncated} (for $n \geq -1$) if the diagonal $f: c \to c \times_d c$ is $(n-1)$-truncated. 
	An object is $n$-truncated if the map to the final object is $n$-truncated. In particular, a morphism $c \to d$ is $(-1)$-truncated (also called mono) if $\Delta:c \to c \times_d c$ is an equivalence. We denote the full sub-category of $n$-truncated objects by $\tau_n\C$ and say $\C$ is an {\it $(n+1,1)$-category} if $\tau_n\C\simeq \C$. 
	\item \label{item:filtered colimit} Let $F: I \to \widehat{\cat}_\infty$ be a diagram, where $I$ is filtered. Then $I$ is weakly contractible \cite[Lemma 5.3.1.18]{lurie2009htt} and so if $F$ is constant its colimit is constant. Moreover, by \cite[Proposition 5.3.3.3]{lurie2009htt} and \cite[Lemma 2.8]{rasekh2021filterquotient}, the colimit commutes with finite limits. As a result, if there is a functor $H:I \times [1] \to \cat_\infty$, with $H(-,0) \simeq F$, $H(-,1) \simeq \D$ constant and $H(i,-)$ mono, then the induced map of colimits $\colim_I F \to \D$ is also mono.  
	\item \label{item:cart rfib} For an $\infty$-category $\C$, a {\it Cartesian fibration} over $\C$ models a presheaf from $\C$ valued in $\infty$-categories and a {\it right fibration} over $\C$  (also called discrete Cartesian fibrations in \cite[Section 5.5]{riehlverity2021elements}) models a presheaf  from $\C$ valued in spaces. They were first introduced in the context of quasicategories  \cite[Section 2.4]{lurie2009htt}, but have been generalized to an arbitrary $\infty$-cosmos \cite[Chapter 5]{riehlverity2021elements}, in a way that coincides with the original definition when we restrict to the $\infty$-cosmos of quasi-categories \cite[Appendix F]{riehlverity2021elements}.
	For a fixed $\infty$-category $\C$, we denote the full sub-quasi-category of $(\widehat{\cat_\infty})_{/\C}$ consisting of Cartesian fibrations by $\Cart_{/\C}$ and the full subcategory of right fibrations   by $\RFib_{/\C}$. Note $\Cart_{/\C}$ (up to categorical equivalence) is independent of the choice of $\infty$-cosmos, as biequivalences of $\infty$-cosmoi preserve and reflect (discrete) Cartesian fibrations \cite[Proposition 10.3.6(x),(xi)]{riehlverity2021elements}.
	\item \label{item:rep right fib} For every $\infty$-category $\C$ and object $c$, there exists a {\it representable right fibration} $\pi_c:\C_{/c} \to \C$ \cite[Corollary 5.5.13]{riehlverity2021elements}. Moreover, we have the {\it Yoneda lemma}: For a right fibration $\R \to \C$ there is an equivalence of spaces $\Map_{(\cat_\infty)_{/\C}}(\C_{/c},\R) \simeq \Fib_c\R$, where $\Fib_c\R$ is the fiber of $\R$ over $c$  \cite[Theorem 5.7.1]{riehlverity2021elements}. In particular, taking $\R= \C_{/d}$, we get an equivalence $\Map_{(\cat_\infty)_{/\C}}(\C_{/c},\C_{/d}) \simeq \Map_\C(c,d)$ \cite[Corollary 5.7.16]{riehlverity2021elements}.
	\item \label{item:cso} Let $\C$ be a finitely complete $\infty$-category. A simplicial object $\X_\bullet:\DD^{op} \to \C$ is a {\it complete Segal object} if it satisfies the {\it Segal condition}: for all $n \geq 2$ the map
	$$\X_n \to \X_1 \times_{\X_0} ... \times_{\X_0} \X_1,$$
	is an equivalence \cite[Definition 3.1]{rasekh2017cartesian}, and the {\it completeness condition}: the square 
	\begin{center}
		\begin{tikzcd}
			\X_0 \arrow[r] \arrow[d] \arrow[dr, phantom, "\ulcorner", very near start]& \X_3 \arrow[d] \\
			\X_1  \arrow[r] & \X_1 \times_{\X_0} \X_1 \times_{X_0} \X_1 
		\end{tikzcd}
	\end{center}
	is a pullback square \cite[Definition 3.3]{rasekh2017cartesian}, where the bottom corner is defined as the limit of the diagram 
	$$\X_1 \xrightarrow{ \ d_0 \ } \X_0 \xleftarrow{ \ d_0 \ } \X_1 \xrightarrow{ \ d_1 \ } \X_0  \xleftarrow{ \ d_1 \ } \X_1.$$
	Complete Segal objects are a model of internal $\infty$-categories and so have a notion of objects, morphisms, equivalences, ... as described in \cite[Section 3]{rasekh2017cartesian}, which we will not require.
	\item \label{item:rep cart fib} For a given complete Segal object $\X$, there is a Cartesian fibration $\C_{/\X_\bullet}$ \cite[Proposition 2.4, Proposition 3.4]{rasekh2017cartesian}, which models the functor $\Map_\C(-,\X_\bullet): \C^{op} \to \widehat{\cat_\infty}$ \cite[Notation 2.5]{rasekh2017cartesian}. We have a {\it Yoneda lemma for representable Cartesian fibrations}: For two complete Segal objects $\X_\bullet,\Y_\bullet$, there is an equivalence $\Map_{\Cart_{/\C}}(\C_{/X_\bullet},\C_{/\Y_\bullet}) \simeq \Map_{\Fun(\DD^{op},\C)}(\X_\bullet,\Y_\bullet)$ \cite[Theorem 2.7]{rasekh2017cartesian}.
	\item \label{item:underlying right fib} There is a fully faithful functor $\RFib_\bullet: \Cart_{/\C} \to \Fun(\DD^{op},\RFib_{/\C})$ with essential image complete Segal objects in right fibrations, where $\RFib_n$ can be explicitly described as the right fibration associated to the functor
	$$\C^{op} \to \widehat{\cat_\infty} \xrightarrow{ \ (-)^{[n]} \ } \widehat{\cat_\infty} \xrightarrow{ \ (-)^\simeq \ } \widehat{\s}.$$
	This is the result of \cite[Section 4]{rasekh2021cartfibcss} and \cite[Section 3]{rasekh2021cartfibmarkedvscso}, and has also been reviewed in \cite[Section 5]{rasekh2021univalence}. The particular case of $\RFib_0$ is also known as the {\it underlying right fibrations} \cite[Corollary 2.4.2.5]{lurie2009htt}.
	\item \label{item:mono Cart fib} A functor of Cartesian fibrations $F:\D \to \E$ over $\C$ is fully faithful if it is fiber-wise fully faithful $\Fib_c\D \hookrightarrow \Fib_c\E$, which, by the argument in the proof of \cite[Proposition 7.6]{rezk2001css} is equivalent to the following 
	being a pullback square of Cartesian fibrations 
	\begin{center}
		\begin{tikzcd}[row sep=0.5in, column sep=1in]
			\RFib_1\D \arrow[r, "\RFib_1F"] \arrow[d] \arrow[dr, phantom, "\ulcorner", very near start] & \RFib_1\E \arrow[d] \\
			\RFib_0\D \times_\C \RFib_0\D \arrow[r, "(\RFib_0F\comma \RFib_0F)"] & \RFib_0\E \times_\C \RFib_0\E
		\end{tikzcd}.
	\end{center}
	\item \label{item:target fib} Let $\C$ be a finitely complete $\infty$-category. Then the target fibration from the arrow $\infty$-category $\O_\C \to \C$ is a Cartesian fibration \cite[Lemma 6.1.1.1]{lurie2009htt}, which represents the functor $\C_{/-}:\C^{op} \to \widehat{\cat_\infty}$, which takes an object $c$ to the over-category $\C_{/c}$. We denote the underlying right fibration $\Und_0\O_\C$ by $\Oall_\C$, and $\Und_1\O_\C$ by $\Aall_\C$ and, by \fibref{item:underlying right fib}, they have value $(\C_{/c})^\simeq$ and $((\C_{/c})^{[1]})^\simeq$, respectively. Moreover, if $S$ is a class of morphisms in $\C$ closed under pullbacks, then we denote the sub-fibration of $\O_\C$ with objects morphisms in $S$, by $\O_\C^S$ and the analogous sub-fibration of $\Oall_\C$ by $\OS_\C$. As a particular case, denote the full subcategory $\O_\C$ ($\Oall_\C$) with objects $n$-truncated morphisms, by $\tau_n\O_\C$ ($\tau_n\Oall_\C$).
	\item \label{item:univalence} Let $\C$ be a finitely complete $\infty$-category. Then a morphism $p:E \to B$ is {\it univalent} if the functor $\C_{/B} \to \Oall_\C$ given by \fibref{item:rep right fib} is fully faithful. Moreover, by \cite[Theorem 4.4]{rasekh2021univalence}, if $\C$ is locally Cartesian closed and $p: E \to B$ is univalent, then there exists a complete Segal object $\N(p): \DD^{op} \to \C$, such that 
	\begin{itemize}
		\item $\N(p)_0\simeq B$
		\item $\N(p)_1 \simeq [E \times B, B \times E]_{B \times B}$, the internal mapping object.
		\item There is a fully faithful inclusion of Cartesian fibrations $\C_{/\N(p)} \hookrightarrow \O_\C$, with essential image morphisms that can be obtained as a pullback of $p$.
	\end{itemize}	
	For more details regarding univalence in finitely complete $\infty$-categories see \cite[Section 2]{rasekh2021univalence}.
	\item \label{item:univalence mono} \cite[Proposition 2.5]{rasekh2021univalence} Let $p: E \to B$ be a univalent morphism. Then in the pullback square 
	\begin{center}
		\begin{tikzcd}
			D \arrow[r] \arrow[d, "q"] \arrow[dr, phantom, "\ulcorner", very near start]& E \arrow[d, "p"] \\
			A \arrow[r, "i"] & B 
		\end{tikzcd}
	\end{center}
	$q$ is univalent if and only if $i$ is mono.
	\item \label{item:univalent overcat} Let $\C$ be a finitely complete $\infty$-category and $p: E \to B$ be univalent. Then for any object $c$ in $\C$, $p\times \id_c$ is univalent in $\C_{/c}$. Indeed, for any morphism $q:y \to x \to c$ we have an equivalence 
	$$\Map_{\O_{\C_{/c}}}(q,p \times \id_c) \simeq  \Map_{\O_\C}(q,p)$$
	and so the desired result follows form \fibref{item:univalence}.
	\item \label{item:idempotent} Let $\Idem$ be the $\infty$-category described in \cite[Definition 4.4.5.2, Remark 4.4.5.3]{lurie2009htt}. An $\infty$-category $\C$ is {\it idempotent complete} if every functor $\Idem \to \C$ has a limit \cite[Corollary 4.4.5.14]{lurie2009htt}. By \cite[Remark 4.4.5.3]{lurie2009htt}, $\Idem$ has a single non-degenerate morphism for each dimension. As as result, if $\C$ is an $n$-truncated, a functor $\Idem \to \C$ corresponds to a functor $\tau_n\Idem \to \C$, which has a limit if and only if $\C$ has finite limits. Hence all finitely complete $n$-truncated $\infty$-categories are idempotent complete.
	\item \label{item:idempotent completion} Let $\widehat{\cat}^\idem$ denote the full subcategory of $\widehat{\cat}$ consisting of idempotent complete $\infty$-categories. By \cite[Proposition 5.1.4.2]{lurie2009htt}, the inclusion $\widehat{\cat}^\idem \hookrightarrow\widehat{\cat}$ has a left adjoint $(-)^{\idem}$ known as the {\it idempotent completion} \cite[Definition 5.1.4.1]{lurie2009htt}.
	\item \label{item:presentability} A quasi-category $\C$ is {\it presentable} ({\it $\omega$-accessible}) if it it satisfies the equivalent conditions given in \cite[Theorem 5.5.1.1]{lurie2009htt} (\cite[Proposition 5.4.2.2]{lurie2009htt}). An $\infty$-category is presentable ($\omega$-accessible) if its underlying quasi-category (\fibref{item:cosmos}) is presentable ($\omega$-accessible). We will simply call it accessible to simplify notation. 
	\item \label{item:compactness} Let $\kappa$ be a regular cardinal and $\C$ an accessible $\infty$-category. Then an object $c$ in $\C$ is $\kappa$-compact if $\Map_\C(c,-): \C \to \s$ commutes with $\kappa$-filtered colimits. Let $\Acc_\infty$ denote the very large $\infty$-category  with objects accessible $\infty$-categories with large set of compact objects and morphism functors that preserve filtered colimits and compact objects. There is an equivalence $(-)^{\cmp}: \Acc \to \widehat{\cat}^{\idem}_\infty$ that takes an accessible $\infty$-category to its full subcategory of compact objects and has inverse $\Ind$, the {\it ind completion} \cite[Proposition 5.4.2.15]{lurie2009htt}, which takes $\C$ to the full subcategory of $\sP(\C) = \Fun(\C^{op},\s)$ consisting of presheaves that can be obtained via filtered colimits \cite[Definition 5.3.5.1]{lurie2009htt}. 
	\item \label{item:compact generated} If $\C$ has finite colimits, then $\Ind(\C)$ is presentable \cite[First paragraph of Subsection 5.5.7]{lurie2009htt} and so has a $(-1)$-truncation functor $\C \to \tau_{-1}\C$ \cite[Proposition 5.5.6.18]{lurie2009htt}.
	\item \label{item:accessible over cat} Let $\C$ be an $\infty$-category and $P: \C^{op} \to \s$. Let $\C_{/P} = \C \times_{\sP(\C)} \sP(\C)_{/P}$. Then the argument in \cite[Corollary 5.1.6.12]{lurie2009htt} implies that the natural functor $\sP(\C_{/P}) \to \sP(\C)_{/P}$ is an equivalence. Moreover, assuming $P$ is obtained from a filtered diagram, we can restrict both sides to presheaves obtained by filtered diagrams and get an equivalence $\Acc(\C_{/P}) \to \Acc(\C)_{/P}$.
	\item \label{item:presentable representale}  Let $\C$ be a presentable $\infty$-category. Then a functor $F:\C^{op} \to \s$ is representable if and only if it takes colimits to limits \cite[Proposition 5.5.2.2]{lurie2009htt}. Moreover, an accessible functor of presentable $\infty$-categories $F: \C \to \D$ is a left adjoint if and only if it preserves colimits \cite[Corollary 5.5.2.9]{lurie2009htt}.
\end{enumerate}

\section{Review of Topos Theory} \label{sec:review}
Before we move on to study the concepts necessary to study elementary $\infty$-toposes, we shortly review the three notions of topos that have already been established: {\it elementary toposes}, {\it Grothendieck toposes} and {\it Grothendieck $\infty$-toposes}. We will start with elementary toposes. 

\begin{defone}\label{def:soc cat}
	Let $\C$ be a finitely complete $1$-category. For a given object $x$, let $\Sub_\C(x)$ be the set of isomorphism classes of subobjects of $x$. Notice this gives us a functor $\Sub_\C:\C^{op}\to \Set$ with functoriality given by pullback. Now a {\it subobject classifier} is an object that represents $\Sub_\C$.
\end{defone}

For more details regarding subobjects and subobject classifiers in categories see \cite[Section I.3]{maclanemoerdijk1994topos}.

\begin{defone} \label{def:et}
	An {\it elementary topos} is a locally Cartesian closed category with subobject classifier. Moreover, a {\it logical functor} is a functor of elementary toposes that preserves finite limits, local Cartesian structure and the subobject classifier. we denote the resulting category by $\Log$.
\end{defone}

For (far more) details regarding elementary toposes see \cite{johnstone2002elephanti,johnstone2002elephantsii,maclanemoerdijk1994topos}. Elementary toposes have a special case, known as Grothendieck toposes. Recall that a category $\G$ is a Grothendieck topos if there exists a small category $\C$ and adjunction 
\begin{center}
	\begin{tikzcd}
		\Fun(\C^{op},\Set)  \arrow[r, shift left=1.8, "L", "\bot"'] \arrow[r, shift right=1.8, "I"', hookleftarrow]  & \G 
	\end{tikzcd},
\end{center}
where $I$ is fully faithful and $L$ preserves finite limits \cite[Section III]{maclanemoerdijk1994topos}. This notion has been generalized to the $\infty$-categorical setting by Lurie \cite[Definition 6.1.0.4]{lurie2009htt}.

\begin{defone} \label{def:groth infty}
	An $\infty$-category $\G$ is {\it presentable} if there exists a small $\infty$-category $\C$ and an adjunction 
	\begin{center}
		\begin{tikzcd}
			\Fun(\C^{op},\s)  \arrow[r, shift left=1.8, "L", "\bot"'] \arrow[r, shift right=1.8, "I"', hookleftarrow] & \G
		\end{tikzcd},
	\end{center}
    such that $L$ is accessible, and is a {\it Grothendieck $\infty$-topos} if $L$ is additionally left exact. 
\end{defone}

\begin{nottwo} 
	Grothendieck $\infty$-toposes are simply called {\it $\infty$-toposes} in \cite{lurie2009htt}, however, we will always use the additional ``Grothendieck" as we want to generalize it in later sections. On the other side, they are called {\it model toposes} in \cite{rezk2010toposes}, as Rezk is using model categorical techniques, and we will not use this terminology as well. 
\end{nottwo}

The definition of Grothendieck $\infty$-topos here is very dependent on the particular choice of $\C$ and we want a more inherent definition, motivated by work of Rezk \cite{rezk2010toposes}, which is more amenable to eventual generalizations. This require us to review some concepts.  Recall that colimits in an $\infty$-category $\G$ are universal if for every morphism $f: x \to y$ the pullback functor $f^*: \G_{/y} \to \G_{/x}$ preserves colimits \cite[Definition 6.1.1.2]{lurie2009htt}. Now we have the following first result.

\begin{propone}
	Colimits in a Grothendieck $\infty$-topos are universal \cite[Theorem 6.1.0.6]{lurie2009htt}, and so, by \fibref{item:presentable representale}, $\G$ is locally Cartesian closed (\fibref{item:lccc}). 
\end{propone}

We want to generalize \cref{def:soc cat} to Grothendieck $\infty$-toposes.

\begin{defone} \label{def:sub}
	Let $\C$ be a finitely complete $\infty$-category. Let $\Sub_\C: \C^{op} \to \Set$ be the composition of $\C_{/-}$ (corresponding to the fibration $\O_\C$ from \fibref{item:target fib}) and the $(-1)$-truncation functor $\tau_{-1}:\cat_\infty \to \Set$ (\fibref{item:truncated}). A subobject classifier is an object $\Omega$ that represents $\Sub_\C$.
\end{defone}

We now have the following in a Grothendieck $\infty$-topos \cite[Subsection 6.1.6]{lurie2009htt}.

\begin{propone} \label{prop:groth subobj}
 $\G$ has a subobject classifier, meaning an object $\Omega: \G^{op} \to \Set$ and a natural isomorphism $\Sub_\C \cong \Hom_\C(-,\Omega)$.
\end{propone} 

We want to move on to the next key property of Grothendieck $\infty$-toposes, {\it local morphisms}. Let $\C$ be a finitely complete $\infty$-category and let $S$ be a set of morphisms in $\C$ closed under pullbacks. Recall from \fibref{item:target fib} that $\O^S_\C$ the full subcategory of $\O_\C$ consisting of objects in $S$ and similarly, $\OS_\C$ the analogous subcategory of $\Oall_\C$. We now have the following definition.

\begin{defone} \label{def:local}
	Let $\C$ be a finitely complete $\infty$-category. Moreover, let $S$ be a set of morphism in $\C$ closed under pullbacks. We say $S$ is {\it local} if the inclusion functor $\OS_\C \to \O^S_\C$ preserves finite colimits.
\end{defone}

Recall that the initial object in $\O_\C$ is just the identity map on the initial object, which is always included in $S$ (as $S$ is closed under pullbacks) hence, the locality condition reduces to $\OS_\C$ being closed under pushouts, meaning for a given pushout square in $\O^S_\C$
\begin{center}
	\begin{tikzcd}
		f \arrow[r, Rightarrow, "\alpha"] \arrow[d, Rightarrow, "\beta"'] & g \arrow[d, Rightarrow, "\gamma"] \\
		h \arrow[r, Rightarrow, "\delta"'] & k \arrow[ul, phantom, "\ulcorner", very near start]
	\end{tikzcd}
\end{center}
if $\alpha$ and $\beta$ are in $\OS_\C$ (meaning they are pullback squares), then $\gamma$ and $\delta$ are also pullback squares. We now have the following result \cite[Theorem 6.1.0.6, Theorem 6.1.3.9]{lurie2009htt}.

\begin{propone}
	Let $\G$ be a locally Cartesian closed presentable $\infty$-category. Then $\G$ is an $\infty$-topos if and only if the set of all morphisms in $\G$ is local.
\end{propone}

We move on to another key property of $\infty$-toposes, {\it descent}. 

\begin{defone} \label{def:descent}
	Let $\C$ be a finitely complete $\infty$-category. Then $\C$ satisfies descent if $(\G_{/-})^\simeq:\G^{op} \to \s$ takes colimits to limits. 
\end{defone}

We now have the following result \cite[Lemma 6.1.3.7, Theorem 6.1.3.9]{lurie2009htt}.

\begin{propone} \label{prop:groth infty topos descent}
	Let $\G$ be a Grothendieck $\infty$-topos, then $\G$ satisfies descent. Moreover, if $\G$ is presentable with universal colimits and satisfies descent, then it is a Grothendieck $\infty$-topos.
\end{propone}

Finally, a key property regarding $\infty$-toposes are {\it universes}. Let $\G$ be an $\infty$-topos and $\kappa$ be a large enough cardinal. Denote by $(\G_{/-})^\kappa$ the full subcategory of $\kappa$-compact objects in $\C_{/-}$ (\fibref{item:compactness}). Then, by \cref{prop:groth infty topos descent}, $((\G_{/-})^\kappa)^\simeq:\G^{op} \to \s$ takes colimits to limits and so, by \fibref{item:presentable representale}, must be represented by an object $\U^\kappa$. This motivates the following definition.

\begin{defone}
	Let $\sP$ be a presentable $\infty$-category and $\kappa$ a cardinal. A {\it universe} or {\it object classifier} for $\kappa$-compact objects is an object $\U^\kappa$ representing the functor $((\sP_{/-})^\kappa)^\simeq:\sP^{op} \to \s$.
\end{defone}

The argument before the definition implies that $\infty$-toposes have universes. However, we do in fact, have the opposite result \cite[Theorem 6.1.6.8]{lurie2009htt}.

\begin{propone}
	A presentable $\infty$-category with universal colimits is a Grothendieck $\infty$-topos if and only if for every large enough cardinal $\kappa$ there exists a universe $\U^\kappa$.
\end{propone}

Up until now, we have seen several different ways of characterizing a Grothendieck $\infty$-topos. We can combine them all into the following theorem.

\begin{theone} \label{the:presentable theorem}
	Let $\G$ be a presentable $\infty$-category. Then the following are equivalent:
	\begin{enumerate}
		\item $\G$ is an $\infty$-topos.
		\item $\G$ is locally Cartesian closed and satisfies descent.
		\item $\G$ is locally Cartesian closed and has object classifiers $\U^\kappa$ for large enough $\kappa$.
		\item $\G$ is locally Cartesian closed and the set of morphisms in $\G$ is local.
	\end{enumerate}
\end{theone}

We can also define an analogue of Grothendieck toposes for $(n,1)$-categories (\fibref{item:truncated}) and concretely we have the following result as given in \cite[Theorem 6.4.1.5]{lurie2009htt}.

\begin{theone} \label{the:presentable theorem n}
	Let $\G$ be a presentable $(n,1)$-category. Then the following are equivalent:
		\begin{enumerate}
		\item $\G$ is an $(n,1)$-topos.
		\item $\G \simeq \tau_{n-1}\hat{\G}$ (\fibref{item:truncated}), where $\hat{\G}$ is a Grothendieck $\infty$-topos.
		\item $\G$ is locally Cartesian closed and has $(n-2)$-truncated object classifiers $\U^\kappa_{\leq n}$ for large enough $\kappa$, meaning an object that represents the functor $\tau_{n-2}((\G_{/-})^\kappa)^\simeq: \G^{op} \to \s_{\leq n-2}$.
		\item $\G$ is locally Cartesian closed and set of $(n-2)$-truncated morphisms in $\G$ is local.
	\end{enumerate}
\end{theone}

\begin{remone} \label{rem:groth subobj n}
	Similar to \cref{prop:groth subobj}, every Grothendieck $(n,1)$-topos $\G$ has a subobject classifier.
\end{remone}
Our goal in the coming sections is to generalize these results from presentable $\infty$-categories to finitely complete $\infty$-categories.

\section{Descent, Cartesian Closures and Universes} \label{sec:universes}
In the last section we reviewed a relation between universes, descent and Cartesian closure in the setting of presentable $\infty$-categories. In this section we want to generalize this relation by introducing a generalization of universes, {\it complete Segal universes}. This will culminate in \cref{the:main univ}, which tells us precisely how these notions interact in the presentable and non-presentable setting.

 As a first step we need to generalize universes from the presentable setting in order to avoid the reliance on set-theoretical assumptions. 

\begin{defone} \label{def:universe}
	Let $\C$ be a finitely complete $\infty$-category. A {\it universe} is a tuple  $(\U,i)$, where $\U$ is an object in $\C$ and $i:\C_{/\U} \hookrightarrow \Oall_\C$ is an inclusion of right fibrations.
\end{defone}

By the Yoneda lemma (\fibref{item:rep right fib}) any functor $i:\C_{/\U} \to \Oall_\C$ is (homotopically) uniquely determined by a choice of morphism $p_\U: \U_* \to \U$, which we call the {\it universal fibration}. Now \fibref{item:univalence} directly implies the following result.

\begin{lemone} \label{lemma:univalent vs universe}
	$(\U,i)$ is a universe if and only if the corresponding universal fibration $p_\U: \U_* \to \U$ is univalent. 
\end{lemone}

We can use this insight to understand the collection of universes.

\begin{defone} \label{def: univ c}
 Let $\C$ be a finitely complete $\infty$-category. Denote by $\Univ_\C$ the (possibly large) poset with objects universes and $\U \leq \U'$ if $\Map_{\Oall_\C}(p_\U,p_{\U'})$ is non-empty, which is equivalent to the existence of a diagram of right fibrations
 \begin{center}
 	\begin{tikzcd}
 		\C_{/\U} \arrow[rr, dashed] \arrow[dr, "i"', hookrightarrow] & & \C_{/\U'} \arrow[dl, "i'", hookrightarrow] \\
 		& \Oall_\C
 	\end{tikzcd}.
 \end{center}
\end{defone} 

We can use this poset to define sufficient universes.

\begin{defone} \label{def:sufficient universes}
	Let $\C$ be a finitely complete $\infty$-category. Then $\C$ has {\it sufficient universes} if the induced map $\ds\coprod_{\U \in \Univ_\C} \C_{/\U} \to \Oall_\C$ is an essentially surjective functor of $\infty$-categories. 
\end{defone}

Concretely, this means that for every morphism $f: Y \to X$ in $\C$, there exists a universe $\U \in \Univ_\C$ and a homotopically unique pullback square 
\begin{equation} \label{eq:classified}
	\pbsq{Y}{\U_*}{X}{\U}{}{f}{p_\U}{\chi_f}
\end{equation}

\begin{notone}
 If such a pullback square exists then we say $\U$ classifies $f$.
\end{notone}

We can combine the definition with \fibref{item:univalence mono} to obtain an elegant classification result for univalent morphisms that generalizes \cite[Corollary  3.10]{gepnerkock2017univalence} from $\infty$-toposes to finitely complete $\infty$-categories with sufficient universes.

\begin{lemone} \label{lemma:universe univalent}
	Let $\C$ be a finitely complete $\infty$-category with sufficient universes. Then a morphism $f$ is univalent if and only if the classifying map $\chi_f$ defined in \ref{eq:classified} is mono.
\end{lemone}

\begin{proof}
	Follows directly from combining the pullback square \ref{eq:classified} with \fibref{item:univalence mono}.
\end{proof}

	Notice the theorem does not depend on the choice of universe $\U$ in $\Univ_\C$ as the collection of universes form a poset. Next we observe that our definition of universe behaves well with respect to over-categories.

\begin{lemone} \label{lemma:universes overcat}
 Let $\C$ be a finitely complete $\infty$-category, let $(\U,i)$ be a universe in $\C$ and let $c$ be an object in $\C$. Then $(\pi_2: \U \times c \to c, i)$ is a universe in $\C_{/c}$. Moreover, if $\C$ has sufficient universes then $\C_{/c}$ has sufficient universes.
\end{lemone} 

\begin{proof}
	First we observe that $\pi_2: \U \times c \to c$ is a universe in $\C_{/c}$. Indeed, this follows directly from \cref{lemma:universe univalent} and the fact that the universal fibration $p_\U \times \id_c: \U_* \times c \to \U \times c$ is univalent in $\C_{/c}$ (\fibref{item:univalent overcat}). Now, let us assume $\C$ has sufficient universes and let $f: d \to e$ be an arbitrary morphism over $c$. By assumption, $f$ is a pullback of $p_\U: \U_* \to \U$, for some $\U$ in $\Univ_\C$ and so we have the pullback square 
	\begin{center}
		\pbsq{d}{\U_* \times c}{e}{\U \times c}{}{}{}{}
	\end{center}
	over $c$, finishing the proof.	
\end{proof}

In certain circumstances we can characterize the existence of sufficient universes in terms that closer resemble a representability condition.

\begin{lemone} \label{lemma:univ filtered}
	Let $\C$ be a finitely complete $\infty$-category with finite coproducts and sufficient universes. Then $\Univ_\C$ is filtered.
\end{lemone}

\begin{proof}
	Let $\U,\U'$ be two universes. As $\C$ has sufficient universes, there exists a universe $\V$ that classifies $p_{\U} \coprod p_{\U'}$, meaning we have pullback diagrams
	\begin{center}
		\begin{tikzcd}[row sep=0.3in, column sep=0.3in]
			\U_* \arrow[r] \arrow[d, "p_\U"]  \arrow[dr, phantom, "\ulcorner", very near start] & \U_* \coprod \U_*' \arrow[r] \arrow[d, "p_\U \coprod p_{\U'}"]  \arrow[dr, phantom, "\ulcorner", very near start] & \V_* \arrow[d, "p_\V"] \\
			\U \arrow[r, hookrightarrow] & \U \coprod \U' \arrow[r, hookrightarrow] & \V 
		\end{tikzcd}.
	\end{center}
	As $p_\U,p_\V$ are univalent the map $\U \to \V$ is mono (\fibref{item:univalence mono}), meaning $\U \leq \V$ in $\Univ_\C$ and similarly $\U' \leq \V$, hence we are done.
\end{proof}

\begin{lemone} \label{lemma:universes representing}
	Let $\C$ be a finitely complete $\infty$-category with finite coproducts. Then $\C$ has sufficient universes if and only if $\ds\underset{\U \in \Univ_\C}{\colim} \C_{/\U} \to \Oall_\C$ is an equivalence of $\infty$-categories.
\end{lemone}

\begin{proof}
	Assume $\C$ has sufficient universes then, by \cref{lemma:univ filtered}, $\Univ_\C$ is filtered. As the diagram is filtered, the colimit is still a full sub-category of $\Oall_\C$ (\fibref{item:filtered colimit}). Moreover, the existence of sufficient universes means it is also essentially surjective, proving it is an equivalence. 
	
	On the other hand, if $\ds\colim_{\U \in \Univ_\C} \C_{/\U} \to \Oall_\C$ is an equivalence, then it is in particular essentially surjective and so for every object $f$ in $\Oall_\C$, there exists a universe $\U$ that classifies $f$, proving $\C$ has sufficient universes.
\end{proof}

General universes help us recover certain objects in the category, however, we often also want to those objects to be closed under certain categorical constructions. 

\begin{defone} \label{def:cp closed}
	Let $\C$ be an $\infty$-category. {\it A universal property} $\cP$ is a collection of functors $P: \C^{op} \to \s$. We say $\C$ is {\it closed under the universal property $\cP$} or {\it satisfies $\cP$} if all functors $P$ are representable in $\C$.
\end{defone}

Examples of such universal properties includes (finite) (co)limits and Cartesian closure.

\begin{defone} \label{def:cp closed univ}
	Let $\C$ be a finitely complete $\infty$-category closed under the universal property $\cP$ and let $\U$ be a universe. We say $\U$ is closed under $\cP$, if for all objects $c$, the full subcategory of $\C_{/c}$ classified by $\U$ is closed under $\cP$. We denote the sub-poset of $\cP$-closed universes of $\C$ by $\Univ_\C^{\cP}$.
\end{defone}

Let us present one example explicitly.

\begin{exone}
	Let $\C$ be a locally Cartesian closed $\infty$-category and a $\U$ a universe. We say $\U$ is Cartesian closed if it closed under the universal property of Cartesian closure. 
\end{exone}

The existence of universes has valuable implications about the $\infty$-category, however, as proven in \cref{lemma:universes representing}, they at best can represent the right fibration $\Oall_\C$ rather than the actual Cartesian fibration $\O_\C$. Hence, in order to obtain a better representability result we need to generalize universes appropriately.

\begin{defone}
	Let $\C$ be a finitely complete $\infty$-category. A {\it complete Segal universe} is a pair $(\U_\bullet: \DD^{op} \to \C, i)$, where $\U_\bullet$ is a simplicial object in $\C$ and $i: \C_{/\U_\bullet} \to \O_\C$ is an inclusion of Cartesian fibrations. Here $\C_{/\U_\bullet}$ is the representable Cartesian fibration defined in \fibref{item:rep cart fib}.
\end{defone}

\begin{warnone}
	Notice our definition of complete Segal universe does not coincide with complete Segal objects internal to the category $\Univ_\C$ (as defined in \fibref{item:cso}). Indeed $\Univ_\C$ is a poset and so any complete Segal object would be constant. 
\end{warnone}

Despite the possible confusion we can justify our naming convention via the following two results that relate complete Segal universes both to universes and complete Segal objects.

\begin{lemone}
	Let $\C$ be a finitely complete $\infty$-category. There is a functor $\Und: \CSUniv_\C \to \Univ_\C$ that takes a complete Segal universe $(\U_\bullet,i)$ to the universe $(\U_0,i)$, which we call the {\it underlying universe}.
\end{lemone}

\begin{proof}
 It suffices to show that $\Und$ is functorial. If $\U_\bullet \leq \U_\bullet'$, then there is an inclusion of Cartesian fibration $\C_{/\U} \to \C_{/\U'}$, which induces an inclusion of right fibrations $\C_{/\U_0} \to \C_{/\U_0'}$, proving that $\U_0 \leq \U_0'$.
\end{proof} 

\begin{lemone}
	Let $(\U_\bullet,i)$ be a complete Segal universe. Then $\U_\bullet$ is a complete Segal object.
\end{lemone}	

\begin{proof}
	By definition $\C_{/\U_\bullet}$ is a Cartesian fibration and so the result follows from \fibref{item:rep cart fib}.
\end{proof}

\begin{defone}
	Let $\C$ be a finitely complete $\infty$-category. Denote by $\CSUniv_\C$ the poset with objects complete Segal universes and $\U_\bullet \leq \U'_\bullet$ if there is a diagram of Cartesian fibrations
	 \begin{center}
		\begin{tikzcd}
			\C_{/\U_\bullet} \arrow[rr, dashed] \arrow[dr, "i"', hookrightarrow] & & \C_{/\U'_\bullet} \arrow[dl, "i'", hookrightarrow] \\
			& \O_\C
		\end{tikzcd}.
	\end{center}
\end{defone} 

\begin{defone}
	Let $\C$ be a finitely complete $\infty$-category. Then $\C$ has {\it sufficient complete Segal universes} if the map $\ds\coprod_{\U_\bullet \in \CSUniv_\C} \C_{/\U_\bullet} \to \O_\C$ is an essentially surjective map of $\infty$-categories. 
\end{defone}

We now have the result analogous to \cref{lemma:universes overcat} for over-categories.

\begin{lemone} \label{lemma:csuniverses overcat}
	 Let $\C$ be a finitely complete $\infty$-category, let $(\U_\bullet,i)$ be a complete Segal universe in $\C$ and let $c$ be an object in $\C$. Then $(\pi_2: \U_\bullet \times c \to c, i)$ is a complete Segal universe in $\C_{/c}$. Moreover, if $\C$ has sufficient complete Segal universes then $\C_{/c}$ has sufficient complete Segal universes.
\end{lemone} 

We also have results analogous to \cref{lemma:univ filtered} and \cref{lemma:universes representing} for complete Segal universes.

\begin{lemone} \label{lemma:csuniv filtered}
	Let $\C$ be a finitely complete $\infty$-category with finite coproducts and sufficient complete Segal universes. Then $\CSUniv_\C$ is filtered.
\end{lemone}

\begin{lemone} \label{lemma:csuniverses representing}
	Let $\C$ be a finitely complete $\infty$-category with finite coproducts. Then $\C$ has sufficient complete Segal universes if and only if $\ds\underset{\U \in \CSUniv_\C}{\colim} \C_{/\U} \to \O_\C$ is an equivalence of $\infty$-categories.
\end{lemone}

We can also additionally assume that complete Segal universes in $\C$ are closed, generalizing \cref{def:cp closed univ}.

\begin{defone} \label{def:csu p closed}
	Let $\C$ be a finitely complete $\infty$-category closed under the universal property $\cP$ and let $\U_\bullet$ be a complete Segal universe. We say $\U_\bullet$ is closed under $\cP$, if for all objects $c$, the full subcategory of $\C_{/c}$ classified by $\U_\bullet$ is closed under $\cP$. We denote the sub-poset of $\cP$-closed complete Segal universes of $\C$ by $\CSUniv_\C^{\cP}$.
\end{defone}

By \cref{def:cp closed univ}, we immediately have the following lemma.

\begin{lemone} \label{lemma:univ cp closed under}
	Let $\C$ be a finitely complete $\infty$-category closed under the universal property $\cP$. Then a complete Segal universe $\U_\bullet$ is closed under $\cP$ if and only if the underlying universe is closed under $\cP$.
\end{lemone}

Moreover, complete Segal universes can have an additional property that we need later on.

\begin{defone} \label{def:csu cart closed}
	Let $\C$ be a finitely complete $\infty$-category and let $\U_\bullet$ be a complete Segal universe. We say $\U_\bullet$ is a {\it locally small} if $\U_0$ classifies the map $(d_1,d_0): \U_1 \to \U_0 \times \U_0$.
\end{defone}

 	The terminology is motivated by locally small categories, which might not have small sets of objects $\cO$ and morphisms $\cM$, but the source target projection $(s,t):\cM \to \cO \times \cO$ is small \cite[Definition 1.1.7]{riehl2017context}.
	Before we can prove the main result, we have the following key proposition.

\begin{propone} \label{prop:lift}
	Let $\C$ be a locally Cartesian closed $\infty$-category. Then the functor $\Und:\CSUniv_\C \to \Univ_\C$ has an inverse $\Lift: \Univ_\C \to \CSUniv_\C$. Moreover, for a given universal property $\cP$, it restricts to an inverse $\Lift: \Univ_\C^\cP \to \CSUniv_\C^\cP$. In particular, there are sufficient universes (closed under $\cP$) if and only if there are sufficient complete Segal universes (closed under $\cP$).
\end{propone}

\begin{proof}
	Let $\U$ be a universe in $\E$ and $p_\U: \U_* \to \U$ the universal fibration. Let $\N(p_\U)$ be the complete Segal object constructed in \fibref{item:univalence}, which satisfies $\N(p_\U)_0 \simeq \U$. This implies that $\Und\circ\Lift$ is the identity. On the other side, for a given complete Segal universe $\U_\bullet$, $\Lift\Und(\U_\bullet)$ is a complete Segal universe with an equivalence of Cartesian fibrations over $\C$
	$$\C_{/\Lift\Und(\U_\bullet)} \simeq \O^\U_\C \simeq \C_{/\U_\bullet} ,$$
	where $\O^\U_\C$ is the full subcategory of $\O_\C$ consisting of morphisms classified by $\U$. Hence, by the Yoneda lemma for complete Segal objects (\fibref{item:rep cart fib}) $\Lift\Und\U_\bullet \simeq \U_\bullet$, which proves that $\Lift\Und$ is also the identity. This proves that $\Und$ and $\Lift$ are inverses.
	
	Now, let $\cP$ be a universal property. Then, by \cref{lemma:univ cp closed under}, a complete Segal universe $\U_\bullet$ is $\cP$-closed if and only if $\Und(\U)$ is $\cP$-closed. Hence, the bijection $(\Und,\Lift)$ between $\Univ_\C$ and $\CSUniv_\C$ restricts to a bijection of posets $\Univ^{\cP}_\C$ and $\CSUniv_\C^{\cP}$.
\end{proof}

We are now ready to state and prove the main theorem.

\begin{theone} \label{the:main univ}
	Let $\E$ be a finitely complete $\infty$-category with finite coproducts closed under the (possibly empty) universal property $\cP$. Consider the following statements:
	\begin{enumerate}
		\item $\E$ is locally Cartesian closed with sufficient universes closed under $\cP$.
		\item $\E$ has sufficient Cartesian closed complete Segal universes closed under $\cP$.
		\item Colimits in $\E$ are universal and the morphisms in $\E$ are local.
		\item Colimits in $\E$ are universal and $\E$ satisfies descent.
	\end{enumerate}
	Then we have the following implications:
	\begin{itemize}
		\item $(1)$ and $(2)$ are equivalent and, assuming they hold, the complete Segal universes are Cartesian closed if and only if they are locally small.
		\item $(3)$ and $(4)$ are equivalent.
		\item $(1),(2)$ imply $(3),(4)$.
		\item If $\E$ is presentable and $\kappa$-compact morphisms are closed under $\cP$ for $\kappa$ large enough, then $(3),(4)$ imply $(1),(2)$.
	\end{itemize}
\end{theone}

\begin{proof}
	$(1) \Leftrightarrow (2)$ If $(1)$ holds, then $(2)$ follows directly from \cref{prop:lift}. On the other side, the existence of sufficient universes closed under $\cP$ follows directly from \cref{lemma:univ cp closed under}. Hence we only need to prove that $\E$ is locally Cartesian closed. Before we can do so we need to construct the appropriate right fibration. 
	
	Let $\Aall_\C \to \C$ be the right fibration introduced in \fibref{item:target fib}. Notice, the source and target projection map gives us a functor of right fibrations $\Aall_\C \to \Oall_\C \times_\C \Oall_\C$. Finally, by the Yoneda lemma, a functor of right fibrations $\C_{/x} \to \Oall_\C \times_\C \Oall_\C$ is given by a choice of two morphisms $f:y \to x,g:z \to x$ in $\C$. The resulting pullback $\Aall_\C \times_{\Oall_\C \times_\C \Oall_\C} \C$ is the right fibration which takes an object $c$ to $\Map_{/c}(y \times_x c,z \times_x c)$ and by \cite[Proposition 2.1]{gepnerkock2017univalence} is representable if and only if the internal mapping object $[y,z]_x$ exists. 
	
	Let $f:y \to x,g:z \to x$ be two arbitrary objects in $\C_{/x}$ and fix a complete Segal universe $\U_\bullet$ that classifies $f \coprod g$ and hence both $f,g$. Denote by $[y,z]_x$ the pullback of $(f,g):x \to \U_0 \times \U_0$ along $(d_1,d_0):\U_1 \to \U_0 \times \U_0$ and notice we now have the following pullback diagram of right fibrations over $\E$
	\begin{center}
		\begin{tikzcd}
			\E_{/[y,z]_x} \arrow[r] \arrow[d] \arrow[dr, phantom, "\ulcorner", very near start] & \E_{/\U_1} \arrow[d] \arrow[r, hookrightarrow] \arrow[dr, phantom, "\ulcorner", very near start] & \Aall_\E \arrow[d] \\
			\E_{/x} \arrow[r, "(f \comma g)"] & \E_{/\U_0} \times_\E \E_{/\U_0} \arrow[r, hookrightarrow] & \Oall_\E \times_\E \Oall_\E 
		\end{tikzcd}
	\end{center}
	where the left hand side is a pullback square by definition of $[y,z]_x$ and the right hand side is a pullback square as the functor $\E_{/\U_\bullet} \to \O_\E$ is fully faithful (\fibref{item:mono Cart fib}). This means the rectangle is a pullback, which proves that $[y,z]_x$ is the internal mapping object of $y,z$ over $x$. This proves that $(2)$ implies $(1)$.
	
	Finally, if $(1),(2)$ hold, then, by \cref{prop:lift}, a complete Segal object $\U_\bullet$ is of the form $\Lift\U$, where $\U$ is a universe and in particular $\U_1 \simeq [\U_* \times \U,\U \times \U_*]_{\U \times \U}$. Hence, if the complete Segal object is Cartesian closed (\cref{def:csu p closed}), then $[\U_* \times \U,\U \times \U_*]_{\U \times \U} \to \U \times \U$ is classified by $\U$, proving $\U_\bullet$ is locally small (\cref{def:csu cart closed}). On the other side, if $\U_1 \to \U_0 \times \U_0$ is classified by $\U_0$, then for all $f:y \to x, g: z \to x$ classified by $\U_0$, $[y,z]_x$ is the pullback of $\U_1 \to \U_0 \times \U_0$ and so also classified by $\U_0$, meaning $\U_\bullet$ is also Cartesian closed. 
	
	\medskip 
	
	$(3) \Leftrightarrow (4)$ This follows directly from \cite[Lemma 6.1.3.7]{lurie2009htt}. Notice the lemma assumes presentability, however, the proof relies on \cite[Lemma 6.1.3.5]{lurie2009htt}, which never uses the presentability assumption.
	
	\medskip 
	
	$(3) \Rightarrow (1)$ Assuming $\E$ is presentable, this follows directly from \cref{the:presentable theorem}.
\end{proof}

\begin{remone}
	In the proof of \cref{the:main univ} if $\E$ has sufficient complete Segal universes, by \cref{lemma:csuniverses overcat}, the over-category $\E_{/x}$ also has sufficient complete Segal universes for all objects $x$. Hence it would have sufficed to prove that $\E$ is Cartesian closed and the case for $\E_{/x}$ would have followed directly. However, our approach above has the additional benefit of giving us an explicit construction of the internal mapping objects as a pullback.
\end{remone}

We will use these various characterizations in the next section to define elementary versions of $\infty$-toposes.

\section{Elementary \texorpdfstring{$\infty$}{oo}-Topos Theory} \label{sec:eht}
 In this section we give a definition of an elementary $\infty$-topos., using the motivation from  \cref{sec:review} and results from \cref{sec:universes}.
 
 \begin{defone} \label{def:eht}
  Let $\cP$ be a universal property. A {\it $\cP$-closed elementary $\infty$-topos} $\E$ is a finitely (co)complete $\infty$-category that satisfies $\cP$ with subobject classifier and sufficient $\cP$-closed complete Segal universes in $\E$. 
 \end{defone} 

 \cref{the:main univ} permits the following alternative characterization as well as valuable implications.
 
 \begin{propone} \label{prop:eht lccc}
 	Let $\cP$ be a universal property. A finitely complete $\infty$-category that satisfies $\cP$ with subobject classifier is a $\cP$-closed elementary $\infty$-topos if and only if it is locally Cartesian closed and has sufficient $\cP$-closed universes. In particular all colimits are universal.
 \end{propone}
 
 \begin{remone} \label{rem:comparison shulman}
 	\cref{prop:eht lccc} implies that the our definition of elementary $\infty$-topos closed under Cartesian closure and finite (co)limits precisely coincides with the definition suggested by Shulman \cite{shulman2017eht}.
 \end{remone}
 
 \begin{nottwo}
 	We will often suppress the notation of the universal property $\cP$ in order to simplify notation.
 \end{nottwo}

 \cref{the:main univ} also has the following valuable implication.
 
 \begin{corone}
 	Let $\E$ be an elementary $\infty$-topos. Then $\E$ satisfies descent (\cref{def:descent}) and the morphisms in $\E$ are local (\cref{def:local}). 
 \end{corone}

Besides universes elementary $\infty$-toposes can also be characterized via univalence, giving us a connection to homotopy type theory \cite[Subsection 2.10]{hottbook2013}.

\begin{corone}
	Let $\cP$ be a universal property. A finitely cocomplete locally Cartesian closed  $\infty$-category that satisfies $\cP$ with subobject classifier $\E$ is an elementary $\infty$-topos if and only if it has sufficient univalent morphisms.
\end{corone}

We can use this result to characterize all univalent morphisms, generalizing \cite[Corollary 3.10]{gepnerkock2017univalence} from Grothendieck $\infty$-toposes, as discussed in \cref{lemma:universe univalent}.

\begin{corone} \label{cor:eht univalent classified}
 Let $\E$ be an elementary $\infty$-topos. A map $f: x \to y$ is univalent if and only if any classifying map $\chi: y \to \U$ given in \ref{eq:classified} is mono.
\end{corone}

\begin{remone} \label{rem:univalent}
	In \cite[Section 6]{rasekh2021univalence} we showed that we can already classify all mono univalent maps in any $\infty$-category with subobject classifier and so in particular elementary toposes. However, we also showed that there are non mono univalent maps in elementary toposes that cannot be classified this way \cite[Example 6.9]{rasekh2021univalence}. The result above shows that an elementary $\infty$-topos is a suitable generalization as it allows us to classify all univalent morphisms.
\end{remone}

We will discuss several examples (and non-examples) of elementary $\infty$-toposes in \cref{sec:examples}, however, in order to help understanding we will give one very explicit example right now.

\begin{exone} \label{ex:space eht}
	Let $\s$ be the $\infty$-category of spaces. We want to show that it satisfies the conditions given in \cref{def:eht}. The existence of finite limits and colimits is evident and the subobject classifier is given by the two element discrete space $\{0,1\}$ and so we focus on showing there are sufficient complete Segal universes. 
	
	Let $\kappa$ be a large enough cardinal. Then, $\s^\kappa$ is a small $\infty$-category and by \fibref{item:css} corresponds to a small complete Segal space, which we simply denote by $\CSS(\s^\kappa)_\bullet: \DD^{op} \to \s$. Notice, the two functors 
	$$\Map(-,\CSS(\s^\kappa)), \CSS(\s_{/-})^\kappa: \s \to \Fun(\DD^{op},\s)$$ 
	are colimit preserving and take the point to $\CSS(\s^\kappa)$. As any two colimit preserving functors out of $\s$ which agree on the point are equivalent \cite[Theorem 5.1.5.6]{lurie2009htt}, we get the desired equivalence of $\infty$-categories $\Map(-,\CSS(\s^\kappa)) \simeq (\s_{/-})^\kappa$. This proves that $\CSS(\s^\kappa)$ is the desired complete Segal universe. Hence, assuming there are sufficiently large cardinals, there are also sufficient complete Segal universes, proving $\s$ is an elementary $\infty$-topos.
\end{exone}

Let us move on to the {\it fundamental theorem of topos theory} \cite{mclarty1992topos} in the $\infty$-categorical setting, also generalizing the analogous result for Grothendieck $\infty$-toposes \cite[Proposition 6.3.5.1]{lurie2009htt}. 

\begin{theone} \label{the:eht fundamental}
	Let $\cP$ be a universal property and let $\E$ be a $\cP$-closed elementary $\infty$-topos and $x$ an object. Then $\E_{/x}$ is also a $\cP$-closed elementary $\infty$-topos. 
\end{theone}

\begin{proof}
	We need to check that $\E_{/x}$ satisfies the conditions in \cref{def:eht}. The case for finite limits and colimits is evident. The case for subobject classifier follows analogous to \cite[Theorem VI.7.1]{maclanemoerdijk1994topos}. The case for complete Segal universes follows from \cref{lemma:csuniverses overcat}.
\end{proof}

We end this section with various implications regarding elementary $\infty$-toposes. First of all we have natural number objects.

\begin{defone}\label{def:nno}
	Let $\C$ be a finitely complete $\infty$-category. A triple $(\bN,s:\bN \to \bN, o:1 \to \bN)$ is a {\it natural number object} if for every triple  $(x,u:x \to x, b:1 \to x)$ in $\C$, the space of of morphisms $f: \bN \to X$ that fit into the diagram 
	\begin{center}
		\begin{tikzcd}[row sep=0.1in, column sep=0.3in] 
			& \bN \arrow[r, "s"] \arrow[dd, "f", dashed] & \bN \arrow[dd, "f", dashed] \\
			1  \arrow[ur, "0"] \arrow[dr, "b"'] \\
			& x \arrow[r, "u"] & x
		\end{tikzcd}
	\end{center}
	is contractible.
\end{defone}

For a more detailed explanation of the definition see \cite[Subsection 2.1]{rasekh2021nno}

\begin{theone}[\citenno{Theorem 4.1.2}] \label{the:eht nno}
	Every elementary $\infty$-topos has a natural number object. 
\end{theone} 

Notice this result is trivial if the elementary $\infty$-topos has countable colimits, as $\coprod_{\bN} 1$ is always a natural number object \cite[Proposition 5.3.2]{rasekh2021nno}. However, there are non-trivial examples, as we shall see in \cref{ex:eht n filter product}. One key implication of the existence of natural number objects is the existence of truncations. To each natural number $n: 1 \to \bN$ we can associate a sphere $S^n$ in $\E$ and a subcategory of $n$-truncated objects $\tau_n\E$ (as the $S^n$-local objects) \cite[Definition 3.5]{rasekh2018truncations}. We now have the following result.

\begin{theone}[\citetrunc{Corollary 4.26}] \label{the:eht ntrunc}
	Let $\E$ be an elementary $\infty$-topos with sufficient Cartesian closed universes closed under finite (co)limits. Let $n: 1 \to \bN$ be a natural number. Then the inclusion functor $\tau_n\E \hookrightarrow \E$ has a left adjoint $\tau_n: \E \to \tau_n\E$.
\end{theone}

We will use this result extensively when studying elementary $(n,1)$-toposes in \cref{sec:eht n}.

\section{The \texorpdfstring{$\infty$}{oo}-Category of Elementary \texorpdfstring{$\infty$}{oo}-Toposes} \label{sec:log infty}
In the last section we defined elementary $\infty$-toposes. We now want to move on and construct a (very large) $\infty$-category of elementary $\infty$-toposes, which requires introducing an appropriate notion of functor. In \cref{def:et} we observed that the appropriate functors of elementary toposes, the {\it logical functors}, preserved finite limits, Cartesian closure and subobject classifier. Taking this as a motivation we want to define $\infty$-logical functors as functors that preserve the appropriate structure. 

Notice that a functor $F: \C_1 \to \C_2$ induces a functor of arrow categories $\sO_F: \sO_{\C_1} \to \sO_{\C_2}$. With this in mind we introduce the following terminology.

\begin{defone} \label{def:infty logical functor}
	Let $\cP$ be a universal property. A finitely (co)continuous functor $F: \E_1 \to \E_2$ of $\cP$-closed elementary $\infty$-toposes is an {\it $\infty$-logical functor} if $\O_F: \O_{\E_1} \to \O_{\E_2}$ restricts to a map of posets $\CSUniv^{\cP}_F: \CSUniv^{\cP}_{\E_1} \to \CSUniv^{\cP}_{\E_2}$, which preserves the subobject classifier.
\end{defone} 

The second condition concretely means that for a given complete Segal universe $\U_\bullet$ in $\CSUniv^{\cP}_{\E_1}$, the simplicial object $F(\U_\bullet)$ is a complete Segal universe in $\CSUniv^{\cP}_{\E_2}$.

 \begin{nottwo}
	Similar to the previous section we will often suppress the notation of the universal property $\cP$ in order to simplify notation.
\end{nottwo}
Recall that logical functors of elementary toposes involves preservation of power objects. As before we can deduce this condition rather than assuming it. Recall that a functor $F: \C_1 \to \C_2$ is locally Cartesian closed if there is a natural equivalence $F[Y,Z]_X \simeq [FY,FZ]_{FX}$.

\begin{propone} \label{prop:logical is lccc}
 Let $F: \E_1 \to \E_2$ be an $\infty$-logical functor of elementary $\infty$-toposes. Then $F$ is locally Cartesian closed. 
\end{propone}

\begin{proof}
	First of all, by \cref{lemma:univ filtered}, $\Univ_{\E_1}$ is filtered. Hence, the equivalence $\colim_{\U \in \Univ_{\E_1}} (\E_1)_{/\U} \simeq \Oall_{\E_1}$ (from \cref{lemma:universes representing}) induces an equivalence
	$$\colim_{\U \in \Univ_{\E_1}} (\E_1)_{/\U \times \U} \simeq \colim_{\U \in \Univ_{\E_1}} \E_{/\U} \times_\E \E_{/\U} \simeq \Oall_{\E_1} \times_{\E_1} \times \Oall_{\E_1}.$$
	Now, by the proof in \cref{the:main univ}, the functor $[-,-]_{-}: \Oall_{\E_1} \times_{\E_1} \times \Oall_{\E_1} \to \E_1$ is given by pulling back along $(d_1,d_0):\U_1 \to \U_0 \times \U_0$. Now the fact that $F$ commutes with pullbacks implies that the following diagram commutes
	\begin{center}
		\begin{tikzcd}[row sep=0.3in, column sep=0.3in]
			\Oall_{\E_1} \underset{\E_1}{\times} \Oall_{\E_1} \arrow[r, "(d_1 \comma d_0)^*"] \arrow[d, "\Oall_F \underset{F}{\times} \Oall_F"'] & \E_1 \arrow[d, "F"] \\
			\Oall_{\E_2} \underset{\E_2}{\times} \Oall_{\E_2} \arrow[r, "(d_1 \comma d_0)^*"] & \E_2
		\end{tikzcd},
	\end{center}
	which proves that $F$ is locally Cartesian closed.
\end{proof}

Using this result we can give an alternative characterization of $\infty$-logical functors, similar to \cref{prop:eht lccc}.

\begin{propone} \label{prop:logical iff lccc}
	Let $\cP$ be a universal property.	A finitely (co)continuous functor $F: \E_1 \to \E_2$ of elementary $\infty$-toposes is an $\infty$-logical functor if and only if $\Oall_F: \Oall_{\E_1} \to \Oall_{\E_2}$ restricts to a map of posets $\Univ^{\cP}_F: \Univ^{\cP}_{\E_1} \to \Univ^{\cP}_{\E_2}$, which preserves the subobject classifier, and $F$ is a locally Cartesian closed functor.
\end{propone}

\begin{proof}
	If $F: \E_1 \to \E_2$ is an $\infty$-logical functor, then by \cref{prop:logical is lccc}, it is a locally Cartesian closed functor. On the other side, we need to prove for every complete Segal universe $\U_\bullet$ in $\E_1$, $F(\U_\bullet)$ is a complete Segal universe. By assumption $F(\U_0)$ is a universe in $\E_2$. Moreover, by \fibref{item:univalence}, $\U_1 \simeq [(\U_0)_* \times \U_0,\U_0\times (\U_0)_*]_{\U_0 \times \U_0}$ and so by assumption $F(\U_1) \simeq [F(\U_0)_* \times F\U_0,F\U_0\times F(\U_0)_*]_{F\U_0 \times F\U_0}$, which, again by \fibref{item:univalence}, proves that $F\U_\bullet$ is a complete Segal universe. This proves that $F$ is an $\infty$-logical functor.
\end{proof}

Similarly, we have a characterization of $\infty$-logical functors via univalent morphisms.

\begin{propone} \label{prop:logical iff univalent}
	A finitely (co)continuous functor $F: \E_1 \to \E_2$ of elementary $\infty$-toposes is an $\infty$-logical functor if and only if $F$ preserves univalent morphisms and the subobject classifier and is a locally Cartesian closed functor.
\end{propone}

\begin{proof}
 By \cref{prop:logical iff lccc} it suffices to prove that if $F$ is (co)continuous, locally Cartesian closed and preserves the subobject classifier , then $F$ preserves universes if and only if it preserves univalent morphisms. However, this follows directly from \cref{lemma:universe univalent}.
\end{proof}

Finally, an $\infty$-logical functor also interacts well with the natural number object.

\begin{propone}
	Let $F: \E_1 \to \E_2$ be $\infty$-logical and $(\bN,s,o)$ a natural number object in $\E_1$. Then $(F\bN,Fs,Fo)$ is a natural number object in $\E_2$.
\end{propone}

\begin{proof}
	By \cite[Definition 2.1.4, Theorem 4.1.2]{rasekh2021nno} a natural number object can be characterized as a triple such that $(s,o): 1 \coprod \bN\to \bN$ is an equivalence and the coequalizer of $s$ and $\id_\bN$ is the terminal object. These two properties are preserved by $F$, as $F$ preserves finite limits and colimits. Hence the result follows.
\end{proof}

Let us give a class of examples of $\infty$-logical functors.

\begin{propone}
	Let $\E$ be an elementary $\infty$-topos and $f: x \to y$ a morphism. Then $f^*: \E_{/y} \to \E_{/x}$ is logical.
\end{propone}

\begin{proof}
	$\E$ preserves finite limits by definition. Moreover, it preserves finite colimits because $\E$ is locally Cartesian closed and hence has a right adjoint. Finally, $f$ preserves complete Segal universes, as $f^*(\U_\bullet \times y \to y) = \U_\bullet \times x \to x$ and \cref{lemma:csuniverses overcat}.
\end{proof}

Notice the identity functor is $\infty$-logical and it is closed under composition giving us the following definition.

\begin{defone}
	Let $\cP$ be a universal property. Let $\Log^{\cP}_\infty$ be the (non-full) subcategory of $\widehat{\cat}_\infty$ with objects $\cP$-closed elementary $\infty$-toposes and morphisms $\infty$-logical functors.
\end{defone}

For later sections it is helpful to have some intermediary notions between general $\infty$-categories and elementary $\infty$-toposes.

\begin{notone} \label{not:infty cats}
	We denote the subcategory of $\widehat{\cat}_\infty$ with objects finitely complete $\infty$-categories with sufficient universes and functors that preserve finite limits and universes by $(\widehat{\cat}_\infty)_\U$. Moreover, denote the subcategory of $(\widehat{\cat}_\infty)_\U$ where objects additionally have a subobject classifier and functors that preserve it by $(\widehat{\cat}_\infty)_{\U,\Prop}$. Finally, denote the subcategory of $(\widehat{\cat}_\infty)_\U$ where the objects are additionally locally Cartesian closed and functors that preserve it by $(\widehat{\cat}_\infty)_{\CSU}$.
\end{notone}

We know the inclusion of Grothendieck $\infty$-toposes and geometric morphisms in the $\infty$-category of $\infty$-categories preserves limits and filtered colimits \cite[Proposition 6.3.2.3, Theorem 6.3.3.1]{lurie2009htt}. Moreover, we have a similar result for the inclusion of elementary toposes in the category of categories (as stated in \cite[Page 218]{maclanemoerdijk1994topos}). We want to observe that $\Log^{\cP}_\infty$ behaves similarly.

\begin{theone} \label{the:filtered colimit eht}
	The inclusion functor $\Log^{\cP}_\infty \to \widehat{\cat}_\infty$ preserves small limits and filtered colimits.
\end{theone}

\begin{proof}
	Let $F:I \to \Log^{\cP}_\infty$ be a diagram. Then the limit in $\widehat{\cat}_\infty$, $\sL =\lim_I F$, has a terminal object given by the tuple of terminal objects. Moreover, for a finite diagram $d:K \to \sL$, by \fibref{item:limits} the limit is given as the terminal object in the pullback $\infty$-category 
	$1 \xrightarrow{ \ d \ } \sL^K \xleftarrow{ \ \Delta \ } \sL,$
	which evidently commutes with limits. The case for finite colimits is analogous. 
	
	Now, for $i \in I$, let $\Omega_i$ be a subobject classifier in $\E_i$ and notice for every morphism $i \to j$ in $I$, the corresponding functor $\E_i \to \E_j$ maps $\Omega_i$ to $\Omega_j$ as the functor is $\infty$-logical. Hence, we have an object $(\Omega_i)_{i \in I}$ in $\sL$. Now, for an arbitrary object $(x_i)_{i \in I}$ we have 
	$$\Map_{\sL}((x_i)_{i \in I},(\Omega_i)_{i \in I}) \simeq \Map_{\E_i}(x_i,\Omega_i)_{i \in I} \cong \Sub_{\E_i}(x_i)_{i \in I} \cong \Sub_{\sL}((x_i)_{i \in I})$$
	which proves that $(\Omega_i)_{i \in I}$ is the subobject classifier. The case for complete Segal universes follows analogously. This proves that $\Log^{\cP}_\infty$ is closed under small limits in $\widehat{\cat}_\infty$.

	Next, we want to observe that the filtered colimit of elementary $\infty$-toposes along $\infty$-logical functors is an elementary $\infty$-topos. Here we use \cite[Theorem 2.26]{rasekh2021filterquotient}. The theorem is stated for filter quotients, however, only uses the fact that the diagram is filtered and that the functors are $\infty$-logical.
\end{proof}

Notice the inclusion $\Log^{\cP}_\infty \to \widehat{\cat}_\infty$ does not preserve the initial object as the initial $\infty$-category is the empty $\infty$-category, which is evidently not an elementary $\infty$-topos. However, we do expect $\Log^{\cP}_\infty$ to have an initial object.

\begin{conjone} \label{conj:free topos}
	$\Log^{\cP}_\infty$ has an initial object.
\end{conjone}

This initial object would be called the {\it free elementary $\infty$-topos}, in analogy to the free elementary topos \cite[Example D4.3.14]{johnstone2002elephantsii}. In fact there is a stronger version of the conjecture, which assigns to every category $\C$ the free elementary $\infty$-topos $L(\C)$.

\begin{conjone} \label{conj:free topos adj}
	For every $\infty$-category $\C$, $(\Log^{\cP}_\infty)_{\C/}$ has an initial object. Equivalently the inclusion $\Log^{\cP}_\infty \to \cat_\infty$ has a left adjoint. 
\end{conjone}

We hope to return to these questions in future work. Finally, note on the other side that the inclusion functor does not preserve the initial object and so cannot have a right adjoint, meaning there cannot be an underlying elementary $\infty$-topos. 

\section{Sufficient vs. External Universes} \label{sec:external universe}
One challenging aspect of the definition of an elementary $\infty$-topos is the existence of sufficient universes as it involves the existence of various classifying objects. It would have been preferable to have a single universe that classifies all objects, however, this would result in size related paradoxes and hence is impossible. 

On the other hand having a better understanding of universes can be of great benefit in particular as we want to eventually strictify the universe and relate it to homotopy type theory (as has been done for Grothendieck $\infty$-toposes by Shulman \cite{shulman2019inftytoposunivalent}). Hence, we here want to introduce a way to obtain such a result using the ind-completion. This requires some aspects of accessible $\infty$-categories. By \fibref{item:compactness}, there is an equivalence $(-)^{\cmp}: \Acc_\infty \to \widehat{\cat}_\infty^{\idem}$, from accessible $\infty$-categories (\fibref{item:presentability}) to idempotent complete $\infty$-categories (\fibref{item:idempotent}). We want to generalize this equivalence to all $\infty$-categories, which requires us to generalize accessible $\infty$-categories.

\begin{defone} \label{def:acc small obj}
	An {\it accessible $\infty$-category with small objects} is a pair $(\A,\C)$, where $\A$ is an accessible category and $\C^{\idem} \simeq \A^{\cmp}$, where $\C^{\idem}$ is the idempotent completion (\fibref{item:idempotent completion}). A {\it functor of accessible $\infty$-categories with small objects}  $F: (\A,\C) \to (\A',\C')$ is an accessible functor $F: \A \to \A'$ that restricts to a functor of small objects $F: \C \to \C'$. We will denote the $\infty$-category of accessible $\infty$-categories with small objects by $\Accsma_\infty$.
\end{defone} 

Notice, there is a fully faithful functor $\Acc_\infty \to \Accsma_\infty$ that takes an accessible $\infty$-category $\A$ to the pair $(\A,\A^{\cmp})$. Similarly $\widehat{\cat}^{\idem}_\infty$ is a full subcategory of $\widehat{\cat}_\infty$. Moreover, let $(-)^{\sma}:\Accsma_\infty \to \widehat{\cat}_\infty$ be the functor that takes $(\A,\C)$ to the category of small objects $\C$ and a functor $F: (\A,\C) \to (\A',\C')$ to its restriction $F: \C \to \C'$. Finally, let $\Ind: \widehat{\cat}_\infty \to \Accsma_\infty$ be the functor that takes $\C$ to $(\Ind\C,\C)$. We now have the following result.

\begin{lemone} \label{lemma:sma ind equiv}
	The functors $(-)^{\sma}: \Accsma_\infty \to \widehat{\cat}_\infty$ and $\Ind:\widehat{\cat}_\infty \to \Accsma_\infty$ are inverses.
\end{lemone}

\begin{proof}
	On the one side $(-)^{\sma}\circ \Ind$ is the identity. On the other side, $\Ind \circ (-)^{\sma}(\A,\C) = (\Ind\C,\C) \simeq (\A,\C)$, where the last step follows from \fibref{item:compactness}.
\end{proof}

This new equivalence does in fact generalize the original equivalence $(\Ind, (-)^{\cmp})$ in the sense that we have the following commutative diagram
	\begin{center}
		\begin{tikzcd}[row sep=0.5in, column sep=0.5in]
			\widehat{\cat}_\infty^{\idem} \arrow[d, hookrightarrow] \arrow[r, shift left=1, "\Ind"] & \Acc_\infty \arrow[d, hookrightarrow] \arrow[l, shift left=1, "(-)^{\cmp}"] \\
			\widehat{\cat}_\infty \arrow[r, shift left=1, "\Ind"] & \Accsma_\infty \arrow[l, shift left=1, "(-)^{\sma}"]
		\end{tikzcd}.
	\end{center}
	
We now use this equivalence to give an alternative characterization of elementary $\infty$-toposes. For a given accessible $\infty$-category with small objects $(\A,\C)$ and object $x$ in $\A$, we use the notation $\C_{/x} \to \C$ for the right fibration obtained via the pullback
\begin{center}
	\begin{tikzcd}
		\C_{/x} \arrow[d, twoheadrightarrow] \arrow[r, hookrightarrow] & \A_{/x} \arrow[d, twoheadrightarrow, "\pi"] \\
		\C \arrow[r, hookrightarrow] & \A 
	\end{tikzcd}.
\end{center}
We want to impose various conditions on the objects and morphisms in $\Accsma$ in order to obtain a subcategory equivalent to $\Log_\infty$.

\begin{defone} \label{def:the universe}
	Let $(\A,\C)$ be an accessible $\infty$-category with small objects. Then it is an {\it accessible $\infty$-category with small objects and (complete Segal) universe}, if $\C$ is closed under finite limits and colimits and there is a (simplicial) object $\U$ ($\U_\bullet$) in $\A$ that represents $\Oall_\C$ ($\O_\C$), meaning there is an equivalence of right (Cartesian) fibrations $\C_{/\U} \simeq \Oall_\C$ ($\C_{/\U_\bullet} \simeq \O_\C$) over $\C$.
\end{defone}
	
We have the following lemma regarding universes.

\begin{lemone} \label{lemma:u universe in a}
	Let $(\A,\C)$ be an accessible $\infty$-category with small objects and universe $\U$. Then the corresponding universal fibration $p_\U: \U_* \to\U$ is univalent in $\A$. Hence $\U$ is a universe in $\A$ as defined in \cref{def:universe}
\end{lemone}

\begin{proof}
	By \fibref{item:accessible over cat}, we have an equivalence $\A_{/\U} \simeq \Ind(\C_{/\U})$. Moreover, by assumption $\C_{/\U} \simeq \Oall_\C$, hence $\A_{/\U} \simeq \Ind(\Oall_\C)$. Now $\Ind(\Oall_\C)$ is the full subcategory of $\Oall_\A$ with objects morphisms $f$ of the form $ \colim_{i\in I} f_i$, where $I \to \Oall_\C$. Hence, $p_\U$ is univalent by \fibref{item:univalence}.
\end{proof}

To better comprehend the role of the universe it is instructive to understand the essential image of $\A_{/\U} \to \Oall_\A$.

\begin{lemone}
	Let $\C$ be closed under finite limits.	An object $f: Y \to X$ in $\Oall_\A$ is in the essential image of $\A_{/\U}$ if and only if for every $Z$ in $\C$ and morphism $Z \to X$, the pullback $Z \times_X Y$ is in $\C$. In particular, a morphism $Y \to 1$ is in the essential image if and only if $Y$ is in $\C$.
\end{lemone}

\begin{proof}
	If $f:Y \to X$ is in $\A_{/\U}$, then $f \simeq \colim_I f_i:Y_i \to X_i$. Let $Z$ be in $\C$, then $Z$ is compact and hence every map $Z \to X$ is of the form $Z \to X_i \to X$ for some $i \in I$. Hence $Z \times_X Y \simeq Z \times_{X_i} Y_i$ which is in $\C$, as $\C$ is closed under finite limits.
	
	On the other side, if $f: Y \to X$ satisfies the condition of the lemma, then 
	$$Y \simeq X \times_X Y \simeq \colim_I( X_i \times_X Y),$$
	where the last step follows from the fact that filtered colimits commute with finite limits (\fibref{item:filtered colimit}). Hence, $f \simeq \colim_I \pi_1:X_i \times_X \times Y \to X_i$.
\end{proof}

\begin{remone}\label{rem:small morphisms}
	We should think of the morphisms classified by $\U$ in $\A$ as ``small morphisms". In \cite{abss2014etuniverse} those are chosen axiomatically as part of the definition of the {\it category of classes}. 
\end{remone}
We can use \cref{lemma:u universe in a} to define functors.

\begin{defone} \label{def:accsma u functor}
	A {\it functor $F:(\A,\C) \to (\A',\C')$ of accessible $\infty$-categories with small objects and universes} is a functor of accessible $\infty$-categories such that $F\U$ is a universe, $F\U \leq \U'$ in $\Univ_{\A'}$ and $F$, if restricted to $\C$, preserves finite limits and colimits. 
	We denote the subcategory of $\Accsma$ with objects accessible $\infty$-categories with small objects and universe and morphism their functors by $(\Accsma_\infty)_\U$. 
\end{defone}

Notice the functors in $(\Accsma_\infty)_\U$ permit following alternative characterization, which follows immediately from \cref{lemma:univalent vs universe}.

\begin{lemone}
	Let $F:(\A,\C) \to (\A',\C')$ be a functor of accessible $\infty$-category with small objects such that $F$ restricted to $\C$ preserves finite limits and colimits. Then $F$ is a functor of accessible $\infty$-categories with small objects and universes if and only if $F\U_*\to F\U$ is univalent and $\U \leq \U'$.
\end{lemone} 

We want to restrict further. For that we need the following lemma. 

\begin{lemone} \label{lemma:acc u vs csu}
	Let $(\A,\C)$ be an accessible $\infty$-category with small objects, such that $\C$ has finite coproducts. Then the following are equivalent:
	\begin{enumerate}
		\item $(\A,\C)$ has a complete Segal universe.
		\item $(\A,\C)$ has a universe and $\C$ is locally Cartesian closed.
		\item $\C$ has sufficient complete Segal universes.
		\item $\C$ has sufficient universes and is locally Cartesian closed.
	\end{enumerate}
\end{lemone}

\begin{proof}
	$(1) \Rightarrow (3)$ If $\C$ has sufficient complete Segal universes, and finite coproducts, then $\CSUniv_\C$ is a directed poset (\cref{lemma:univ filtered}) and so it has a colimit $\U_\bullet = \colim_{(\U_\bullet)_i \in \CSUniv_\C} (\U_\bullet)_i$ and we have 
	$$\C_{/\U_\bullet} \simeq \underset{(\U_\bullet)_i \in \CSUniv_\C}{\colim} \C_{/(\U_\bullet)_i} \simeq \O_\C,$$
	where the last step follows from \cref{lemma:csuniverses representing}, proving that $\U_\bullet$ is a complete Segal universe in $(\Ind\C,\C)$.
	
	\medskip 
	
	$(3) \Rightarrow (1)$ If $(\Ind\C,\C)$ has a complete Segal universe $\U_\bullet$, then by definition there exists an equivalence $\colim_I (\U_\bullet)_i \simeq \U_\bullet$, where $I$ is filtered. By construction of the filtered colimit the map $(\U_\bullet)_i \to \U_\bullet$ is mono and so, by \fibref{item:univalence mono} and \cref{lemma:univalent vs universe}, $(\U_\bullet)_i$ is a complete Segal universe. Now, let $f: y \to x$ be an arbitrary morphism in $\C$, then there exists $x \to \U_0$ that classifies $f$. The fact that $x$ is compact implies that there exists a factorization $x \to (\U_i)_0 \to \U_0$ (\fibref{item:compactness}). Hence, the complete Segal universe $(\U_\bullet)_i$ classifies $f$, proving $\C$ has sufficient complete Segal universes. 
	
	\medskip 
	
	$(2) \Leftrightarrow (4)$ Follows from the same argument as in the previous two steps.
	 
	 \medskip 
	 
	 $(3) \Leftrightarrow (4)$ Follows from \cref{the:main univ}.
\end{proof}

This now motivates the following definition. 

\begin{defone}
	An {\it accessible elementary model} $(\A,\E)$ is an accessible $\infty$-category with small objects that satisfies the equivalent conditions of \cref{lemma:acc u vs csu}.
\end{defone}

We now move on to the appropriate notion of functors.

\begin{lemone} \label{lemma:acc u vs csu functor}
	Let $(\A,\E)$ and $(\A',\E')$ be two accessible elementary models. Let $F: (\A,\E) \to (\A',\E')$ be a functor of accessible $\infty$-categories with small objects such that $F: \E \to \E'$ preserves finite limits and colimits. Then the following are equivalent: 
	\begin{enumerate}
		\item $F$ takes the complete Segal universe $\U_\bullet$ to the complete Segal universe $F(\U_\bullet) \leq \U_\bullet'$.
		\item $F$ takes the universe $\U$ to the universe $F\U \leq \U'$ and $F$ is locally Cartesian closed.
		\item $F: \E \to \E'$ induced a map of complete Segal universes $\CSUniv_F: \CSUniv_{\E} \to \CSUniv_{\E'}$.
		\item $F: \E \to \E'$ induced a map of universes $\Univ_F: \Univ_{\E} \to \Univ_{\E'}$ and $F$ is locally Cartesian closed.
	\end{enumerate}
\end{lemone}

\begin{proof}
	$(1) \Rightarrow (3)$ Let $(\U_\bullet)_i $ be a complete Segal universe in $\E$. By definition of the complete Segal universe $\U_\bullet$ in $\A$, there is a map $(\U_\bullet)_i \to \U_\bullet$, which by \fibref{item:univalence mono}, it mono. As $F$ preserves finite limits it preserves monos and so we have 
	 $$F((\U_\bullet)_i) \hookrightarrow F(\U_\bullet) \hookrightarrow \U'_\bullet,$$
	 which, again by \fibref{item:univalence mono}, proves that $F$ induces a functor $\CSUniv_F: \CSUniv_{\E} \to \CSUniv_{\E'}$.
	 
	 \medskip 
	 
	 $(3) \Rightarrow (1)$ Let $\U_\bullet$ be the complete Segal universe in $\A$ and notice $\U_\bullet \simeq \colim_I (\U_\bullet)_i$, where $I$ is filtered. As $F$ is accessible it preserves filtered colimits and so $F(\U_\bullet) \simeq \colim_I F((\U_\bullet)_i)$. As we have $F((\U_\bullet)_i) \leq \U_\bullet'$, we have $F(\U_\bullet) \leq \U_\bullet'$ by \fibref{item:filtered colimit}, finishing the proof.
	
	\medskip
	
	$(2) \Leftrightarrow (4)$ This follow from an analogous argument to the previous two steps.
	
	\medskip 
	
	$(3) \Leftrightarrow (4)$ Follows from \cref{prop:logical iff lccc}.
\end{proof}

We can now give the desired definition.

\begin{defone} \label{def:accelem infty}
  A {\it functor $F:(\A,\C) \to (\A',\C')$ of accessible elementary models} is a functor of accessible $\infty$-categories that satisfies the equivalent conditions of \cref{lemma:acc u vs csu functor}. We denote the subcategory of $\Accsma_\infty$ with objects accessible elementary models and their functors by $\AccElem_\infty$.
\end{defone}

\begin{defone}
		An accessible $\infty$-category with small objects $(\A,\C)$ has a {\it subobject classifier} if there exists an object $\Omega$ that represents the functor $\Sub_\C$, as defined in \cref{def:sub}, meaning there is a natural isomorphism $\Hom_\A(c,\Omega) \cong \Sub_\C(c)$ for all objects $c$ in $\C$.
\end{defone}

\begin{lemone} \label{lemma:accessible neg truncation}
		Let $(\A,\E)$ be an accessible $\infty$-category with small locally Cartesian closed small objects that has a universe $\U$. Then the inclusion functor $\tau_{-1}\E \to \E$ has a left adjoint $\tau_{-1}:\E \to \tau_{-1}\E$.
\end{lemone}

\begin{proof}
	$\E$ has finite limits and colimits, is locally Cartesian closed, and satisfies descent, where the last condition follows from \cref{the:main univ}. Hence, by \cite[Theorem 4.1.2]{rasekh2021nno}, $\E$ has a natural number object and so in particular a $(-1)$-truncation \cite[Theorem 2.29]{rasekh2018truncations}. 
\end{proof}

\begin{lemone} \label{lemma:subobject classifier}
	Let $(\A,\E)$ be an accessible $\infty$-category with small objects that has a universe $\U$. Then $(\A,\E)$ has a subobject classifier. 
\end{lemone}

\begin{proof}
 By \cref{lemma:accessible neg truncation} there is a functor $\tau_{-1}: \Oall_\E \to \Oall_\E$ over $\E$. We can extend it to a functor $\Ind\tau_{-1}: \Ind\Oall_\E \to \Ind\Oall_\E$. Combining \cref{def:the universe} and \fibref{item:accessible over cat}, we have $\Ind\Oall_\E \simeq \Ind\E_{/\U} \simeq \A_{/\U}$ and so by the Yoneda lemma for right fibrations (\fibref{item:rep right fib}) we get a morphism $\tau_{-1}: \U \to \U$ in $\A$. Now, $\C$ has finite limits and colimits and so, by \fibref{item:compact generated},  $\Ind\C \simeq \A$ is presentable and has a $(-1)$-truncation functor $\tau_{-1}: \A \to \tau_{-1}\A$. Applying it to $\tau_{-1}: \U\to \U$ gives us a factorization $\U \hookrightarrow \Omega \to \U$, such that $\Map(X,\Omega) \simeq \tau_{-1}(\E_{/X}) \cong \Sub_\E(X)$, proving $\Omega$ is a subobject classifier.  
\end{proof}

\begin{defone} \label{def:small subobject}
	An accessible $\infty$-category with small objects $(\A,\C)$ with universe $\U$ satisfies {\it propositional resizing} if the object $\Omega$ is small. Let $(\Accsma_\infty)_{\U,\Prop}$ be the full subcategory of $(\Accsma_\infty)_\U$ consisting of objects that satisfy propositional resizing and functors in $(\Accsma_\infty)_{\U}$ that preserve the small subobject classifier. Define $(\AccElem_\infty)_{\Prop}$ similarly. 
\end{defone}

We are now ready to state and prove the main theorem.

\begin{theone} \label{the:main external stuff}
	There is a diagram of $\infty$-categories
	\begin{center}
		\begin{tikzcd}[column sep=-0.1in]
			& (\widehat{\cat}_\infty)_\U \arrow[rrr, "\simeq" very near start, "\Ind" near end] & & & (\Accsma_\infty)_\U & \\
			(\widehat{\cat}_\infty)_{\CSU} \arrow[ur, hookrightarrow] & & & \AccElem_\infty \arrow[ur, hookrightarrow]& & \\
			& & (\widehat{\cat}_\infty)_{\U,\Prop} \arrow[rrr, "\simeq" very near start, "\Ind" near end] \arrow[uul, hookrightarrow] & & & (\Accsma_\infty)_{\U,\Prop} \arrow[uul, hookrightarrow]\\
			& \Log_\infty \arrow[rrr, "\simeq" near start, "\Ind" near end] \arrow[uul, hookrightarrow] \arrow[ur, hookrightarrow] & & & (\AccElem_\infty)_{\Prop}  \arrow[ur, hookrightarrow] & 
			 \arrow[from=2-1, to=2-4, crossing over, "\simeq" near start, "\Ind" near end]
			 \arrow[from=4-5, to=2-4, hookrightarrow, crossing over]
		\end{tikzcd},
	\end{center}
	where the vertical functors are faithful and the horizontal functors are equivalences. Here the $\infty$-categories on the left hand side were defined in \cref{not:infty cats}.
\end{theone}

\begin{proof}
	By \cref{lemma:sma ind equiv} we have an equivalence $\Ind: \widehat{\cat}_\infty \to \Accsma_\infty$. We need to prove the equivalence restricts appropriately and for that it suffices to prove that $\Ind$ and its inverse $(-)^{\sma}$ preserve the additional structure. 
	
	\begin{enumerate}
		\item {\bf (Co)limits and Cartesian Closure:} First of all if $\C$ has finite (co)limits or is locally Cartesian closed, then the accessible $\infty$-category with small objects $(\Ind\C,\C)$ also has the same property for small objects and the same applies to functors.
		
		\item {\bf Sufficient Universes vs. Universe:} By \cref{lemma:acc u vs csu}, $(\A,\C)$ has a universe if and only if $\C$ has sufficient universes. Moreover, by the analogous argument to the one in \cref{lemma:acc u vs csu functor} (where we ignore the additional locally Cartesian closed assumption), a functor $F: (\A,\E) \to (\A',\E')$ preserves the universe if and only if $F:\E \to \E'$ restricts to a map $\Univ_F: \Univ_{\E_1} \to \Univ_{\E_2}$. 
		
		\item {\bf Subobject classifier vs. Propositional Resizing:} By \cref{def:small subobject}, $\C$ has a subobject classifier if and only if $(\A,\C)$ has a subobject classifier that satisfies propositional resizing. Moreover, it is evident that $F: \E \to \E'$ preserves the subobject classifier if and only if $F: (\A,\E) \to (\A',\E')$ preserves the small subobject classifier, as it is the same object.
	\end{enumerate}
	As we have confirmed that all properties are preserved appropriately, our proof is finished. 
\end{proof}

\begin{remone}
	The equivalence relates elementary $\infty$-toposes with universes and elementary models. One additional question is whether we can give a similar result for elementary $\infty$-toposes with $\cP$-closed universes. This would require using a different approach as the one given in the proof of \cref{the:main external stuff}. Indeed, the proof uses the fact that the universe $\U$ in $(\A,\E)$ can be recovered as a filtered colimit of universes $\U_i$ in $\E$, however, we cannot choose those universes to be closed under $\cP$.
\end{remone}

\begin{remone} \label{rem:justification}
	The definition of $\infty$-logical functor given in \cref{def:infty logical functor} involves an intricate condition because of the existence of sufficient universes. On the other hand, the definition of functors of accessible elementary models is straightforward and very much modeled on the definition of logical functors of elementary toposes. The fact that these two notions coincide (as proven in \cref{the:main external stuff}) gives a strong indication that the definition of $\infty$-logical functor is in fact appropriately chosen.
\end{remone}

\section{Elementary \texorpdfstring{$(n,1)$}{(n,1)}-Topos Theory} \label{sec:eht n}
Up until now our work has focused on generalizing elementary $1$-toposes (\cref{def:et}) to elementary $\infty$-toposes in a way that appropriately generalizes Grothendieck $\infty$-toposes. However, there similarly are Grothendieck $(n,1)$-toposes for all $0 \leq n \leq \infty$, which sit in between the two boundary cases and have been studied extensively in \cite[Section 6.4]{lurie2009htt}. In this section we generalize the results from \cref{sec:eht}, \cref{sec:log infty}  and \cref{sec:external universe}, from $\infty$-categories to $(n,1)$-categories (\fibref{item:truncated}), in particular introducing {\it elementary $(n,1)$-topos theory}.

\begin{nottwo}
	It is an unfortunate fact of higher categorical terminology that an $\infty$-category by now refers to an $(\infty,1)$-category, whereas an $n$-category usually refers to an $(n,n)$-category creating possible confusion. In order to avoid such confusion, we will consistently use the notation $(n,1)$-category, when $n < \infty$. 
\end{nottwo}

We want to start with an appropriate notion of universe for $(n,1)$-categories. As every object in an $(n,1)$-category is $(n-1)$-truncated, for a given morphism $p:E \to B$, the space $\Map_\C(-,B)$ is $(n-1)$-truncated and so the induced map (\fibref{item:rep right fib}) $\C_{/B} \to \Oall_\C$ can only be fully faithful if it takes value in the full subcategory $\tau_{n-2}\Oall_\C$ (\fibref{item:target fib}). This suggests the following notion of universe (\cref{def:universe}) in the setting of $(n,1)$-categories.

\begin{defone} \label{def:n trunc universe}
	Let $\C$ be a finitely complete $(n,1)$-category and $\cP$ a universal property. A $\cP$-closed {\it $(n-2)$-truncated universe} is a pair $(\U,i)$, where $\U$ is an object in $\C$ and $i: \C_{/\U} \to \tau_{n-2} \Oall_\C$ is a fully faithful functor, such that the full subcategory of $\tau_{n-2} \Oall_\C$ with objects in the essential image of $i$ satisfy $\cP$ (\cref{def:cp closed}). 
\end{defone}

Similar to \cref{def: univ c} we denote the poset of $\cP$-closed $(n-2)$-truncated universes by $\Univ_\C^{\cP,n-2}$ and, following \cref{def:sufficient universes}, a finitely complete $(n,1)$-category has sufficient universes if every $(n-2)$-truncated morphism in $\C$ is classified by a universe. We now have the following appropriate result corresponding to \cref{the:main univ}.

\begin{corone} \label{cor:eht n lccc}
	Let $\E$ be a finitely complete locally Cartesian closed $(n,1)$-category with finite coproducts closed under the (possibly empty) universal property $\cP$. If $\E$ has sufficient $(n-2)$-truncated universes closed under $\cP$ then colimits in $\E$ are universal and $(n-2)$-truncated morphisms in $\E$ are local. The opposite holds if $\E$ is presentable and $\kappa$-compact morphisms are closed under $\cP$ for $\kappa$ large enough.
\end{corone}

\begin{proof}
	If $\E$ has sufficient universes closed under $\cP$, then the proof follows using the same argument as in \cref{the:main univ}. On the other side, assuming $\E$ is presentable the other direction follows from \cite[Theorem 6.4.1.5]{lurie2009htt}, as reviewed in \cref{the:presentable theorem n}.
\end{proof}

We now have all the necessary ingredients to define elementary $(n,1)$-toposes.

\begin{defone} \label{def:eht n}
	Let $\cP$ be a universal property. and $n>0$. A finitely (co)complete locally Cartesian closed $(n,1)$-category is an {\it elementary $(n,1)$-topos} if it has a subobject classifier and sufficient $\cP$-closed $(n-2)$-truncated universes.
\end{defone}

This definition is in fact an appropriate generalization of the previous cases. 

\begin{exone}
	If $n=1$, then the subobject classifier itself is precisely a universe for all $(-1)$-truncated morphisms. Hence, the existence of sufficient universes becomes vacuous and elementary $(1,1)$-toposes in the sense of \cref{def:eht n} directly correspond to elementary toposes, as defined in \cref{def:et}.
\end{exone}

\begin{exone}
	If $n= \infty$, then it is a finitely (co)complete locally Cartesian closed $\infty$-category with subobject classifier and sufficient universes for $(n-2)$-truncated morphisms, which corresponds to all morphisms as $n=\infty$ and hence recovers the definition of an elementary $\infty$-topos, by \cref{prop:eht lccc}.
\end{exone}

\begin{remone} \label{rem:eht not global}
	Unlike in the case of $\infty$-toposes, we cannot recover the Cartesian closure from the universes, as morphisms that are not $(n-2)$-truncated are not classified.
\end{remone}

We also have an analogous result to \cref{cor:eht univalent classified}.

\begin{corone} \label{cor:univalent n}
	Let $\E$ be an elementary $(n,1)$-topos. An $(n-2)$-truncated map $f: x \to y$ is univalent if and only if any classifying map $\chi_f: y \to \U$ is mono.
\end{corone}

Notice for $n=1$ \cref{cor:univalent n} recovers the already well-known fact that univalent monos are classified by the subobject classifier, as discussed in \cref{rem:univalent}. There is also a fundamental theorem of elementary $(n,1)$-toposes, with the same proof as in \cref{the:eht fundamental}.
\begin{corone}
	Let $\cP$ be a universal property and let $\E$ be a $\cP$-closed elementary $(n,1)$-topos and $x$ an object. Then $\E_{/x}$ is also a $\cP$-closed elementary $(n,1)$-topos. 
\end{corone}

We move on to the $(n,1)$-categorical analogue of \cref{sec:log infty} and study {\it $(n,1)$-logical functors}.

\begin{defone} \label{def:logical n}
	A finitely (co)continuous locally Cartesian closed functor $F: \E_1 \to \E_2$ of elementary $(n,1)$-toposes is an {\it $(n,1)$-logical functor} if $\O_F: \O_{\E_1} \to \O_{\E_2}$ restricts to a map of posets $\Univ^{\cP,n-2}_F: \Univ^{\cP,n-2}_{\E_1} \to \Univ^{\cP,n-2}_{\E_2}$, which preserves the subobject classifier.
\end{defone} 

\begin{exone}
 In the case $n=1$, this precisely reduces to the case of logical functors of elementary toposes \cref{def:et}, whereas in the case of $n=\infty$ it coincides with $\infty$-logical functors, by \cref{prop:logical iff lccc}.
\end{exone}

Similarly we denote the $\infty$-category with objects $\cP$-closed elementary $(n,1)$-toposes and morphisms $(n,1)$-logical functors by $\Log^{\cP}_{(n,1)}$. We now have the analogous result to \cref{the:filtered colimit eht} with the same proof.

\begin{corone} \label{cor:filtered colimit eht n}
	The inclusion functor $\Log^{\cP}_{(n,1)} \to \widehat{\cat}_{(n,1)}$ preserves small limits and filtered colimits. 
\end{corone}

Finally, we want to adapt the results in \cref{sec:external universe}. First of all, let $(\widehat{\cat}_{(n,1)})_\U,(\widehat{\cat}_{(n,1)})_{\U,\LCCC},(\widehat{\cat}_{(n,1)})_{\U,\Prop}$ be defined analogous to $(\widehat{\cat}_{\infty})_\U,(\widehat{\cat}_{\infty})_{\U,\LCCC},(\widehat{\cat}_{\infty})_{\U,\Prop}$ in \cref{not:infty cats}, respectively. 

On the other hand, following \fibref{item:idempotent}, if $\E$ is a finitely complete $(n,1)$-category, where $n<\infty$, then $\E$ is idempotent complete. Hence, we do not need to use the notion of accessible $\infty$-category with small objects (\cref{def:acc small obj}). Rather we have the following.

\begin{defone}
	Let $\Acc_{(n,1)}$ be the full subcategory of $\Acc_\infty$, consisting of accessible $\infty$-categories such that $\A^{\cmp}$ is an $(n,1)$-category and notice the equivalence in \fibref{item:compactness} restricts to an equivalence $\Ind: \widehat{\cat}_{(n,1)} \to \Acc_{(n,1)}$.
\end{defone}

Let $(\Acc_{(n,1)})_\U,\AccElem_{(n,1)}, \ (\Acc_{(n,1)})_{\U,\Prop}, \ (\AccElem_{(n,1)})_{\Prop}$ be defined as the full subcategories of $(\Accsma_{\infty})_\U$ (\cref{def:accsma u functor}), $\AccElem_{\infty}$ (\cref{def:accelem infty}), $(\Accsma_{\infty})_{\U,\Prop}$, $(\AccElem_{\infty})_{\Prop}$ (\cref{def:small subobject}), respectively, with objects in $\Acc_{(n,1)}$. Now, we have the following result.

\begin{corone}
	Let $n<\infty$. There is a diagram of $\infty$-categories
	\begin{center}
		\begin{tikzcd}[column sep=-0.1in]
			& (\widehat{\cat}_{(n,1)})_\U \arrow[rrr, "\simeq" very near start, "\Ind" near end] & & & (\Acc_{(n,1)})_\U & \\
			(\widehat{\cat}_{(n,1)})_{\CSU} \arrow[ur, hookrightarrow] & & & \AccElem_{(n,1)} \arrow[ur, hookrightarrow]& & \\
			& & (\widehat{\cat}_{(n,1)})_{\U,\Prop} \arrow[rrr, "\simeq" very near start, "\Ind" near end] \arrow[uul, hookrightarrow] & & & (\Acc_{(n,1)})_{\U,\Prop} \arrow[uul, hookrightarrow]\\
			& \Log_{(n,1)} \arrow[rrr, "\simeq" near start, "\Ind" near end] \arrow[uul, hookrightarrow] \arrow[ur, hookrightarrow] & & & (\AccElem_{(n,1)})_{\Prop}  \arrow[ur, hookrightarrow] & 
			\arrow[from=2-1, to=2-4, crossing over, "\simeq" near start, "\Ind" near end]
			\arrow[from=4-5, to=2-4, hookrightarrow, crossing over]
		\end{tikzcd},
	\end{center}
	where the vertical functors are faithful and the horizontal functors are equivalences. 
\end{corone}

We end this section with a discussion how the various truncation levels compare to each other. It is well established that for all $m \leq n \leq \infty$, every $(m,1)$-category is an $(n,1)$-category. Now that we have established a notion of elementary $(n,1)$-topos, we want to understand whether the various dimensions are related similarly. 

\begin{propone} \label{prop:eht m not n}
	Let $-1 \leq n < m \leq \infty$ and $\E$ an $(n,1)$-category. Then $\E$ is not an elementary $(m,1)$-topos.
\end{propone}

\begin{proof}
	By \cref{def:n trunc universe}, a universe for $(m-2)$-truncated objects is necessarily $(m-1)$-truncated, but not $(m-2)$-truncated. Hence, if $n<m$, then every object in $\E$ is $(n-1)$-truncated and so also $(m-2)$-truncated and so cannot be a universe, meaning $\E$ cannot be an elementary $(m,1)$-topos. 
\end{proof}

Hence, elementary $(n,1)$-toposes are sensitive to dimension. However, these notions are still related.

\begin{theone} \label{the:truncation eht}
	Let $\cP$ be the universal property consisting of finite limits, finite colimits and local Cartesian closure. Then for all $n\geq 0$ we have a functor $\tau_n: \Log^{\cP}_\infty \to \Log^{\cP}_{(n+1,1)}$.
\end{theone}

\begin{proof}
  We need to prove that for every $\cP$-closed elementary $\infty$-topos $\E$, $\tau_n\E$ is a $\cP$-closed elementary $(n+1,1)$-topos and for every $\infty$-logical functor $F: \E_1 \to \E_2$, $\tau_nF: \tau_n\E_1 \to \tau_n\E_2$ is an $(n+1,1)$-logical functor. By \cref{the:eht ntrunc}, there is an $n$-truncation functor $\tau_n:\E \to \tau_n\E$, which is left adjoint to the inclusion $\tau_n\E \hookrightarrow \E$. This implies that $\tau_n(\E)$ is closed under finite limits and colimits. Moreover, by \cite[Proposition 2.11]{rasekh2021univalence}, $\tau_n$ is a locally Cartesian closed localization, which implies that $\tau_n\E$ is locally Cartesian closed as well. 
  
  Next, the subobject classifier $\Omega$ is $0$-truncated and so in particular an object in $\tau_n(\E)$. Finally, by \fibref{item:rep right fib}, the truncation functor $\tau_{n-1}: \Oall_\E \to \Oall_\E$ induces a map of universes $\tau_{n-1}: \U\to \U$ and applying epi-mono factorization (which exists by \cref{the:eht ntrunc} for $n=-1$) gives us a factorization $\U \to \U_{\leq n-1} \to \U$. By definition, $\Map_\E(x,\U_{\leq n-1}) \simeq \tau_{n-1}(\E_{/x})$, which proves that $\U_{\leq n-1}$ is a universe for $(n-1)$-truncated objects in $\tau_n\E$. As $\E$ has sufficient universes and this process constructs a $(n-1)$-truncated universe out of every universe, this proves that $\tau_n\E$ also has sufficient $(n-1)$-truncated universes. Hence $\tau_n\E$ is an elementary $(n+1,1)$-topos.
  
  Now, let $F: \E_1 \to \E_2$ be an $\infty$-logical functor. Then $\tau_nF: \tau_n\E_1 \to \tau_n\E_2$ preserves finite limits and colimits and, again by  \cite[Proposition 2.11]{rasekh2021univalence}, the local Cartesian closure. Finally, as $F(\U)$ is a universe and $F$ respects truncation levels, $F(\U_{\leq n-1})$ is again a universe for $(n-1)$-truncated objects, proving $\tau_nF$ is $(n+1,1)$-logical.
\end{proof}

\begin{remone}
	The result from \cref{the:truncation eht} would also hold if $\cP$ is a universal property that includes finite limits, finite colimits and local Cartesian closure and such that for all $n \geq 0$, if $\C$ satisfies $\cP$, then $\tau_n\C$ also satisfies $\cP$.
\end{remone}

This result has some fascinating implications and possibilities for generalizations. By \cref{the:eht ntrunc}, for a given elementary $\infty$-topos $\E$, with natural number $n: 1 \to \bN$ we have an $n$-truncation $\tau_n: \E \to \tau_n\E$.
We would like to consider $\tau_n\E$ an example of an elementary $(n+1,1)$-topos. However, this is not possible if the natural number is non-standard. This suggest the following interesting conjecture.

\begin{conjone}
	Let $\E$ be an elementary $\infty$-topos. Let $\cat_\E$ denote the $\infty$-category of $\E$-enriched $\infty$-categories (as defined by Gepner and Haugseng in \cite{gepnerhaugseng2015enrichedinftycat}) and notice that $\E$ is itself $\E$ enriched \cite[Corollary 7.4.10]{gepnerhaugseng2015enrichedinftycat}. Fix an internal natural number $n$ in the natural number object $\bN$ of $\E$. There is a notion of elementary $(n+1,1)$-topos and applying the $n$-truncation functor $\tau_n:\cat_\E \to \cat_{\E_{\leq n}}$  gives us an example thereof.
\end{conjone}

Notice \cref{the:truncation eht} has an inverse in the context of Grothendieck $\infty$-toposes, meaning every Grothendieck $(n+1,1)$-topos is equivalent to the full subcategory of $n$-truncated objects in a Grothendieck $\infty$-topos, as explained in \cref{the:presentable theorem n}. At this level of generality, this result does not generalize to elementary $(n,1)$-toposes.

\begin{exone}
	The category of finite sets is an elementary topos, but is not equivalent to the category of $0$-truncated objects of an elementary $\infty$-topos. Indeed, by \cref{the:eht nno}, every elementary $\infty$-topos has a natural number object, which is by definition $0$-truncated and so an element in the sub-category of $0$-truncated objects.
\end{exone}

What prevents the existence of a lift is the lack of universes, which by definition exist if $n>1$, suggesting the following conjecture. 

\begin{conjone} \label{conj:truncated eht}
	Let $n>1$. For every elementary $(n+1,1)$-topos $\E$, there exists an elementary $\infty$-topos $\hat{\E}$, such that $\tau_n\hat{\E} \simeq \E$.
\end{conjone}

\section{Examples and Non-Examples} \label{sec:examples}
In this section we cover various examples and non-examples of elementary $(n,1)$-toposes.

\medskip

{\bf Grothendieck $(n,1)$-Toposes:} Let us start with one key example, for which we have the following result.

\begin{propone} \label{prop:grothendieck eht}
	Let $n \leq \infty$, let $\cP$ be a universal property and $\sP$ a presentable $(n,1)$-category that satisfies $\cP$. If there exists sufficient cardinals $\kappa$, such that the full subcategory of $\kappa$-compact objects is closed under $\cP$, then $\sP$ is an elementary $(n,1)$-topos with $\cP$-closed universes if and only if it is a Grothendieck $(n,1)$-topos. 
\end{propone}

\begin{proof}
	If $\sP$ is a Grothendieck $(n,1)$-topos, then, by \cref{the:presentable theorem n}, it is locally Cartesian closed and has sufficient universes $\U^\kappa$ for $\kappa$-small $(n-2)$-truncated objects and a subobject classifier (\cref{rem:groth subobj n}), and hence an elementary $(n,1)$-topos, by \cref{cor:eht n lccc}. On the other side, if $\sP$ is an elementary $(n,1)$-topos, then it is locally Cartesian closed and has sufficient universes for $(n-2)$-truncated objects and so a Grothendieck $(n,1)$-topos, by \cref{the:presentable theorem n}.
\end{proof}

\begin{exone}
	Assuming sufficiently large cardinals, $\s_{\leq n-1}$ is an elementary $(n,1)$-topos. 
\end{exone}

\medskip 

{\bf Finite Spaces:} As the category of finite sets is an elementary topos, we might expect that finite spaces are an elementary $\infty$-topos. However, that is not the case.

\begin{exone} \label{ex:omega compact spaces}
	The category of $\omega$-compact spaces is not an elementary $\infty$-topos. In fact it is not even closed under finite limits and colimits. Indeed, if it has finite limits it has the terminal object $*$. Moreover, if it has finite colimits it has the suspension of $* \coprod*$, which is $S^1$, then finite limits imply it has $\Omega S^1 \cong \mathbb{Z}$, which is certainly not $\omega$-compact.
\end{exone}

\medskip

{\bf $1$-Inaccessible Spaces:} There is a way to correct the previous example that is due to Lo Monaco \cite{monaco2021eht}.  Recall that a cardinal $\kappa_{acc}$ is {\it $1$-inaccessible} if it is weakly inaccessible and for every $\kappa < \kappa_{acc}$, there exists a weakly inaccessible cardinal $\kappa < \kappa_1 < \kappa_{acc}$. We now have the following result as stated in \cite[Proposition 7.3]{monaco2021eht}.

\begin{exone}
	Let $\kappa_{acc}$ be a $1$-inaccessible cardinal. The $\infty$-category of $\kappa_{acc}$-small spaces $\s^{\kappa_{acc}}$ is an elementary $\infty$-topos, using the same argument as in \cref{ex:space eht}. Indeed, for every weakly inaccessible cardinal $\kappa < \kappa_{acc}$, we have a complete Segal universe $\CSS(\s^\kappa)$ (\fibref{item:css}).
\end{exone}

We can in fact generalize the result by Lo Monaco in the expected manner.

\begin{exone}
	Let $\kappa_{acc}$ be a $1$-inaccessible cardinal. The $(n+1,1)$-category of $\kappa_{acc}$-small $n$-truncated spaces $(\s_{\leq n})^{\kappa_{acc}}$ is an elementary $(n+1,1)$-topos. Indeed, as $\kappa_{acc}$ is weakly inaccessible, it is closed under finite limits, finite colimits and locally Cartesian closed. Moreover, $\{0,1\}$ is evidently $\kappa_{acc}$-small and so $(\s_{\leq n})^{\kappa_{acc}}$ has a subobject classifier. Finally, for every weakly inaccessible cardinal $\kappa< \kappa_1$, let $((\s_{\leq n-1})_{\kappa})^\simeq$, be the maximal subgroupoid (\fibref{item:underlyinggroupoid}) of the $\infty$-category of $\kappa$-small $n-2$-truncated spaces. Then by the same argument as the one given in \cref{ex:space eht} we have an equivalence
	$$\Map(X,((\s_{\leq n-1})^{\kappa})^\simeq) \simeq (((\s_{\leq n-1})_{/X})^{\kappa})^\simeq, $$
	which proves that $((\s_{\leq n-1})^{\kappa})^\simeq$ is a universe for $(n-1)$-truncated spaces.
\end{exone} 

Notice, the inclusion functor $(\s_{\leq n})^{\kappa_{acc}} \to \s_{\leq n}$ is the identity functor on all objects, immediately giving us the following result.

\begin{corone}
	The inclusion functor $(\s_{\leq n})^{\kappa_{acc}} \to \s_{\leq n}$ is $(n+1,1)$-logical.
\end{corone}

We can think of this result as a generalization of the fact that the inclusion from finite sets to sets is logical.

\medskip 

{\bf $\pi$-Finite Spaces:} There is an alternative notion of smallness of spaces, known as {\it $\pi$-finite spaces} \cite[Definition E.0.7.8]{lurie2018sag} or {\it truncated coherent spaces} \cite{anel2021eht}, which are spaces for which the disjoint union of all homotopy groups is finite. We have the following result similar to \cref{ex:omega compact spaces}.

\begin{exone}
	The $\infty$-category of bounded coherent spaces does not form an elementary $\infty$-topos, as it again not closed under finite limits and colimits, using the fact that $\mathbb{Z}$ is not bounded coherent.
\end{exone}

\begin{remone}
	The title of \cite{anel2021eht} might create the impression that $\pi$-finite spaces are an elementary $\infty$-topos, however, as the abstract explains, the goal is not to prove that it satisfies all conditions in \cref{def:eht}, but rather that it satisfies many of its conditions (such as finite limits and local Cartesian closure).
\end{remone}

\medskip 

{\bf Kan Objects:} In \cite{ghss2021effective}, the authors construct an $\infty$-category of {\it Kan objects} out of a wide range of categories with countable coproducts and finite limits. While the resulting $\infty$-category always has finite limits and colimits and satisfies descent \cite[Proposition 10.1]{ghss2021effective}, we do in fact have an example that does not have a subobject classifier \cite[Example 11.8]{ghss2021effective}. However, it is not known whether they have sufficient complete Segal universes.

\medskip 

{\bf Filter Quotient $(n,1)$-Toposes:} One interesting example of elementary toposes that are not Grothendieck are filter quotients \cite{adelmanjohnstone1982serreclasses}. They have been generalized to non-trivial examples of elementary $\infty$-toposes in \cite{rasekh2021filterquotient}. We can now use our results to generalize it to elementary $(n,1)$-toposes that are not Grothendieck $(n,1)$-toposes, proving that the opposite to \cref{prop:grothendieck eht} does not hold.

First we need to review several concepts.

\begin{defone} \label{def:filter}
	A filter is a partially ordered set $(\Phi,\leq)$ that satisfies the following three conditions:
	 \begin{enumerate}
		\item $\Phi \neq \emptyset$.
		\item $\Phi$ is downward directed, meaning that for any two object $x, y \in \Phi$ there exists $z \in \Phi$ such that $z \leq x$ and $z \leq y$.
		\item $\Phi$ is upward closed, meaning that if $x \leq y$ and $x \in \Phi$, then $y \in \Phi$.
	\end{enumerate}
	Notice $\Phi^{op}$ is a filtered category and if $\Phi$ has a maximal object $1$, then $(1)$ and $(3)$ implies that $1 \in \Phi$.
\end{defone}

\begin{defone} 
	A finitely complete $\infty$-category with product-closed filter of subobjects is a finitely complete $\infty$-category $\C$ along with a filter $\Phi$ on the poset $\Sub_\C(1)$ that is closed under products. 
\end{defone}

For a given filter, let $\C_{/-}: \Phi^{op} \to \cat_\infty$ be the restriction of $\C_{/-}$ (\fibref{item:target fib}) with the inclusion $\Phi^{op} \to \C^{op}$. Concretely, it is the functor that takes $U \in \Phi$ to $\C_{/\U}$ and $U \leq V$ to $- \times U:\C_{/V} \to \C_{/U}$. We now have the following definition.

\begin{defone}
	Let $(\C,\Phi)$ be a finitely complete $\infty$-category with product-closed filter of subobjects. Define the {\it filter quotient} $\C_{\Phi}$ as the colimit $\colim_{U \in \Phi^{op}}\C_{/U}$ and denote $P_\Phi: \C \simeq \C_{/1} \to \C_{\Phi}$ as the {\it quotient functor}.
\end{defone}

We now have the following result generalizing \cite{rasekh2021filterquotient}. 

\begin{propone}
	Let $(\E,\Phi)$ be a finitely complete $(n,1)$-category with product-closed filter of subobjects. Then $P_\Phi: \E \to \E_\Phi$ 
	preserves 
	\begin{enumerate}
		\item finite limits and colimits
		\item subobject classifiers
		\item locally Cartesian closure
		\item $(n-2)$-truncated universes
	\end{enumerate}
	In particular if $\E$ is an elementary $(n,1)$-topos, then $\E_\Phi$ is one as well and $P_\Phi$ is $(n,1)$-logical
\end{propone}

\begin{proof}
	The diagram $\Phi^{op}$ is filtered (\cref{def:filter}) and so $\E_\Phi$ is a filtered colimit and so the result follows from the same argument given in \cref{cor:filtered colimit eht n}.
\end{proof}

Let us give one explicit example motivated by \cite[Subsection 3.2]{rasekh2021filterquotient}.

\begin{exone} \label{ex:eht n filter product}
	Let $\bN$ be the set of natural numbers and $\Phi$ the Fr{\'e}chet filter of cofinite subsets \cite[Example 1.47]{rasekh2021filterquotient}. Then this induces a filter of subobjects on the $(n+1,1)$-category $(\s_{\leq n})^\bN$ and we call the resulting filter quotient the {\it filter product} elementary $(n+1,1)$-topos and denote it by $\prod_\Phi\s_{\leq n}$. Following the explanation in \cite[Subsection 3.2]{rasekh2021filterquotient}, it does not have countable products and coproducts and the natural number object is non-standard, giving us a very explicit example of an elementary $(n+1,1)$-topos that is not Grothendieck.
\end{exone}

\bibliographystyle{alpha}
\bibliography{main}

\end{document}